

\magnification=1200


\catcode`\@=11
\font\eightrm=cmr8
\font\eighti=cmmi8
\font\eightsy=cmsy8
\font\eightbf=cmbx8
\font\eighttt=cmtt8
\font\eightit=cmti8
\font\eightsl=cmsl8
\font\sixrm=cmr6
\font\sixi=cmmi6
\font\sixsy=cmsy6
\font\sixbf=cmbx6

\font\tensans=cmss10
\font\fivesans=cmss10 at 5pt
\font\sevensans=cmss10 at 7pt
\newfam\sansfam
\textfont\sansfam=\tensans\scriptfont\sansfam=\sevensans
\scriptscriptfont\sansfam=\fivesans
\def\sans{\fam\sansfam\tensans}

\font\petcap=cmcsc10
\font\bftwelve=cmbx10 at 12pt

\def\eightpoint{%
  \textfont0=\eightrm \scriptfont0=\sixrm \scriptscriptfont0=\fiverm
  \def\rm{\fam\z@\eightrm}%
  \textfont1=\eighti \scriptfont1=\sixi \scriptscriptfont1=\fivei
  \def\oldstyle{\fam\@ne\eighti}%
  \textfont2=\eightsy \scriptfont2=\sixsy \scriptscriptfont2=\fivesy
  \textfont\itfam=\eightit
  \def\it{\fam\itfam\eightit}%
  \textfont\slfam=\eightsl
  \def\sl{\fam\slfam\eightsl}%
  \textfont\bffam=\eightbf \scriptfont\bffam=\sixbf
  \scriptscriptfont\bffam=\fivebf
  \def\bf{\fam\bffam\eightbf}%
  \textfont\ttfam=\eighttt
  \def\tt{\fam\ttfam\eighttt}%
  \abovedisplayskip=9pt plus 2pt minus 6pt
  \abovedisplayshortskip=0pt plus 2pt
  \belowdisplayskip=9pt plus 2pt minus 6pt
  \belowdisplayshortskip=5pt plus 2pt minus 3pt
  \smallskipamount=2pt plus 1pt minus 1pt
  \medskipamount=4pt plus 2pt minus 1pt
  \bigskipamount=9pt plus 3pt minus 3pt
  \normalbaselineskip=9pt
  \setbox\strutbox=\hbox{\vrule height7pt depth2pt width0pt}%
  \let\bigf@ntpc=\eightrm \let\smallf@ntpc=\sixrm
  \normalbaselines\rm}

\def\tenpoint{%
  \textfont0=\tenrm \scriptfont0=\sevenrm \scriptscriptfont0=\fiverm
  \def\rm{\fam\z@\tenrm}%
  \textfont1=\teni \scriptfont1=\seveni \scriptscriptfont1=\fivei
  \def\oldstyle{\fam\@ne\teni}%
  \textfont2=\tensy \scriptfont2=\sevensy \scriptscriptfont2=\fivesy
  \textfont\itfam=\tenit
  \def\it{\fam\itfam\tenit}%
  \textfont\slfam=\tensl
  \def\sl{\fam\slfam\tensl}%
  \textfont\bffam=\tenbf \scriptfont\bffam=\sevenbf
  \scriptscriptfont\bffam=\fivebf
  \def\bf{\fam\bffam\tenbf}%
  \textfont\ttfam=\tentt
  \def\tt{\fam\ttfam\tentt}%
  \abovedisplayskip=6pt plus 2pt minus 6pt
  \abovedisplayshortskip=0pt plus 3pt
  \belowdisplayskip=6pt plus 2pt minus 6pt
  \belowdisplayshortskip=7pt plus 3pt minus 4pt
  \smallskipamount=3pt plus 1pt minus 1pt
  \medskipamount=6pt plus 2pt minus 2pt
  \bigskipamount=12pt plus 4pt minus 4pt
  \normalbaselineskip=12pt
  \setbox\strutbox=\hbox{\vrule height8.5pt depth3.5pt width0pt}%
  \let\bigf@ntpc=\tenrm \let\smallf@ntpc=\sevenrm
  \normalbaselines\rm}

\catcode`\@=12

\def\ssm{\mathbin{
\kern1pt\raise0.36ex\hbox{$\scriptscriptstyle\rm\backslash$}}}
\def\restriction{{\kern0.3pt\raise1pt\hbox{$\scriptstyle|$\kern-2.5pt
\raise-2pt\hbox{\rm\`{}}\kern-1.7pt}}}
\def\codim{\mathop{\rm codim}\nolimits}

\def\finpr{\hfill \hbox{
\vrule height 1.453ex  width 0.093ex  depth 0ex
\vrule height 1.5ex  width 1.3ex  depth -1.407ex\kern-0.1ex
\vrule height 1.453ex  width 0.093ex  depth 0ex\kern-1.35ex
\vrule height 0.093ex  width 1.3ex  depth 0ex}}

\def\pointir{\discretionary{.}{}{.\kern.35em---\kern.7em}}
\def\titre#1|{\par\vskip .5cm\penalty -100
              \vbox{\centerline{\bf #1}
                    \vskip 5pt}
              \penalty 500}
\def\bigtitle#1|#2|{\par\vskip .5cm\penalty -100
               \vbox{\centerline{\bf #1}\centerline{\bf #2}
                     \vskip 5pt}
              \penalty 500}
\long\def\tha#1|#2\fintha{\par\vskip 5pt
                    {\petcap #1\pointir}{\sl #2}\par\vskip 5pt}
\def\remarque#1|{\par\vskip5pt{{\sl #1}\pointir}}
\newif\ifchiffre
\def\chiffre{\chiffretrue}
\chiffre
\newdimen\laenge
\def\lettre#1|{\setbox3=\hbox{#1}\laenge=\wd3\advance\laenge by 3mm
\chiffrefalse}
\def\article#1|#2|#3|#4|#5|#6|#7|%
    {{\ifchiffre\leftskip=7mm\noindent
     \hangindent=2mm\hangafter=1
\llap{[#1]\hskip1.35em}{\petcap #2}\pointir {\sl #3}, {\rm #4},
\nobreak{\bf #5} ({\oldstyle #6}), \nobreak #7.\par\else\noindent
\advance\laenge by 4mm \hangindent=\laenge\advance\laenge by -4mm\hangafter=1
\rlap{[#1]}\hskip\laenge{\petcap #2}\pointir
{\sl #3}, #4, {\bf #5} ({\oldstyle #6}), #7.\par\fi}}
\def\livre#1|#2|#3|#4|#5|%
    {{\ifchiffre\leftskip=7mm\noindent
    \hangindent=2mm\hangafter=1
\llap{[#1]\hskip1.35em}{\petcap #2}\pointir{\sl #3}, #4, {\oldstyle #5}.\par
\else\noindent
\advance\laenge by 4mm \hangindent=\laenge\advance\laenge by -4mm
\hangafter=1
\rlap{[#1]}\hskip\laenge{\petcap #2}\pointir
{\sl #3}, #4, {\oldstyle #5}.\par\fi}}
\def\divers#1|#2|#3|#4|%
    {{\ifchiffre\leftskip=7mm\noindent
    \hangindent=2mm\hangafter=1
     \llap{[#1]\hskip1.35em}{\petcap #2}\pointir{\sl #3}, {\rm #4}.\par
\else\noindent
\advance\laenge by 4mm \hangindent=\laenge\advance\laenge by -4mm
\hangafter=1
\rlap{[#1]}\hskip\laenge{\petcap #2}\pointir{\sl #3}, {\rm #4}.\par\fi}}

\def\buildo#1\over#2{\mathrel{\mathop{\null#2}\limits^{#1}}}
\def\buildu#1\under#2{\mathrel{\mathop{\null#2}\limits_{#1}}}

\def\nb{{\rm I\!N}}

\def\pb{{\rm I\!P}}
\def\rb{{\rm I\!R}}

\def\cb{{\mathchoice {\setbox0=\hbox{$\displaystyle\rm C$}\hbox{\hbox
to0pt{\kern0.4\wd0\vrule height0.9\ht0\hss}\box0}}
{\setbox0=\hbox{$\textstyle\rm C$}\hbox{\hbox
to0pt{\kern0.4\wd0\vrule height0.9\ht0\hss}\box0}}
{\setbox0=\hbox{$\scriptstyle\rm C$}\hbox{\hbox
to0pt{\kern0.4\wd0\vrule height0.9\ht0\hss}\box0}}
{\setbox0=\hbox{$\scriptscriptstyle\rm C$}\hbox{\hbox
to0pt{\kern0.4\wd0\vrule height0.9\ht0\hss}\box0}}}}
\def\qb{{\mathchoice {\setbox0=\hbox{$\displaystyle\rm
Q$}\hbox{\raise
0.15\ht0\hbox to0pt{\kern0.4\wd0\vrule height0.8\ht0\hss}\box0}}
{\setbox0=\hbox{$\textstyle\rm Q$}\hbox{\raise
0.15\ht0\hbox to0pt{\kern0.4\wd0\vrule height0.8\ht0\hss}\box0}}
{\setbox0=\hbox{$\scriptstyle\rm Q$}\hbox{\raise
0.15\ht0\hbox to0pt{\kern0.4\wd0\vrule height0.7\ht0\hss}\box0}}
{\setbox0=\hbox{$\scriptscriptstyle\rm Q$}\hbox{\raise
0.15\ht0\hbox to0pt{\kern0.4\wd0\vrule height0.7\ht0\hss}\box0}}}}
\def\tb{{\mathchoice {\setbox0=\hbox{$\displaystyle\rm
T$}\hbox{\hbox to0pt{\kern0.3\wd0\vrule height0.9\ht0\hss}\box0}}
{\setbox0=\hbox{$\textstyle\rm T$}\hbox{\hbox
to0pt{\kern0.3\wd0\vrule height0.9\ht0\hss}\box0}}
{\setbox0=\hbox{$\scriptstyle\rm T$}\hbox{\hbox
to0pt{\kern0.3\wd0\vrule height0.9\ht0\hss}\box0}}
{\setbox0=\hbox{$\scriptscriptstyle\rm T$}\hbox{\hbox
to0pt{\kern0.3\wd0\vrule height0.9\ht0\hss}\box0}}}}
\def\zb{{\mathchoice {\hbox{$\sans\textstyle Z\kern-0.4em Z$}}
{\hbox{$\sans\textstyle Z\kern-0.4em Z$}}
{\hbox{$\sans\scriptstyle Z\kern-0.3em Z$}}
{\hbox{$\sans\scriptscriptstyle Z\kern-0.2em Z$}}}}

\let\la=\longrightarrow
\let\ov=\overline
\let\wt=\widetilde
\let\wh=\widehat

\def\than{\tha\noindent}
\def\rem#1|{\remarque\noindent{\petcap#1}|}
\def\proof{\noindent{\sl{Proof\pointir}\ }}
\def\ld{,\,\ldots\,,}

\def\Hom{{\rm Hom}}

\def\Supp{{\rm Supp}}

\def\Id{{\rm Id}}
\def\ch{{\rm ch}}
\def\pr{{\rm pr}}
\def\Im{{\rm Im}}

\def\op{\overline\partial}
\def\fc{{\cal F}}
\def\ic{{\cal I}}
\def\oc{{\cal O}}
\def\compact{\subset\!\subset}

\def\stimes{\mathop{\kern0.7pt
\vrule height 0.4pt depth 0pt width 5pt\kern-5pt
\vrule height 5.4pt depth -5pt width 5pt\kern-5pt
\vrule height 5.4pt depth 0pt width 0.4pt\kern4.6pt
\vrule height 5.4pt depth 0pt width 0.4pt\kern-6.5pt
\raise0.3pt\hbox{$\times$}\kern-0.7pt}\nolimits}
\def\lhra{\lhook\joinrel\longrightarrow}
\def\srelbar{\mathrel{\smash{\hbox{\rm-}}}}
\def\merto{\srelbar\,\srelbar\,\rightarrow}

\tenpoint
\pretolerance=500 \tolerance=1000 \brokenpenalty=5000
\hoffset=0.66cm 
\hsize=12.5cm
\vsize=19cm
\parskip 5pt plus 1pt
\parindent=0.833cm

\ \vskip36pt
\centerline{\bftwelve Algebraic approximations of holomorphic maps}
\vskip4pt
\centerline{\bftwelve from Stein domains to projective manifolds}

\vskip12pt
{\noindent\hangindent0.6cm\hangafter-1
{\bf Jean-Pierre Demailly$\,$\footnote{${}^1$}{\eightpoint\rm
Research partially supported by Institut Universitaire de France}
\hfill L\'aszl\'o Lempert$\,$\footnote{${}^2$}{\eightpoint\rm
Research partially supported by National Science Foundation Grant No.\
DMS-9303479}\kern0.6cm

}}

{\noindent\hangindent0.6cm\hangafter-4{\sl
Universit\'e de Grenoble I\hfill
Purdue University\kern0.6cm\break
Institut Fourier, BP 74\hfill
Department of Mathematics\kern0.6cm\break
U.R.A. 188 du C.N.R.S.\hfill
West Lafayette, IN 47907, U.S.A.\kern0.6cm\break
38402 Saint-Martin d'H\`eres, France

}}

\vskip8pt
\centerline{\bf Bernard
Shiffman$\,$\footnote{${}^3$}{\eightpoint\rm Research
partially supported by National Science Foundation Grant No.\ DMS-9204037}}
\vskip3pt\centerline{\sl Johns Hopkins University}
\centerline{\sl Department of Mathematics} \centerline{\sl Baltimore, MD
21218, U.S.A.}

\vskip24pt


\vskip5pt
\noindent
{\bf A.M.S. Classification 1985:} 32E30, 32H20, 14C30

\vskip24pt
\centerline{\bf Table of contents}\vskip3pt
{\eightpoint\noindent
1. Introduction\dotfill p.\ 2\break
2. Holomorphic and Nash algebraic retractions\dotfill p.\ 8\break
3. Nash algebraic approximation on Runge domains in affine algebraic
varieties\dotfill p.\ 13\break
4. Nash algebraic approximations omitting ample divisors\dotfill
p.\ 19\break
5. Exhaustion of Stein manifolds by Runge
domains of affine algebraic manifolds\dotfill p.\ 27\break
References\dotfill p.\ 30}

\vfill\break

\titre 1. Introduction|
The present work, which was motivated by the study of the Kobayashi
pseudo\-distance on algebraic manifolds, proceeds from the general
philosophy that analytic objects can be approximated by algebraic objects
under suitable restrictions. Such questions have been extensively studied in
the case of holomorphic functions of several complex variables and can be
traced back to the Oka-Weil approximation theorem (see [We35] and [Oka37]).
The main approximation result of this work (Theorem~1.1) is used to show that
both the Kobayashi pseudodistance and the Kobayashi-Royden infinitesimal
metric on a quasi-projective algebraic manifold $Z$ are computable solely in
terms of the closed algebraic curves in $Z$ (Corollaries~1.3 and~1.4).

Our general goal is to show that algebraic approximation is always
possible in the cases of holomorphic maps to
quasi-projective manifolds (Theorems~1.1 and~4.1) and of locally free sheaves
(Theorem~1.8 and Proposition~3.2). Of course, to avoid difficulties arising
from the possible occurrence of infinite topology or geometry near the
boundary, it is necessary to shrink the objects and to restrict to relatively
compact domains. Since we deal with algebraic approximation, a central notion
is that of Runge domain: By definition, an open set $\Omega$ in a Stein space
$Y$ is said to be a {\sl Runge domain} if $\Omega$ is Stein and if the
restriction map $\oc(Y)\to\oc(\Omega)$ has dense range. It is well known that
$\Omega$ is a Runge domain in $Y$ if and only if the holomorphic hull with
respect to $\oc(Y)$ of any compact subset $K\subset\Omega$ is contained
in~$\Omega$. If~$Y$ is an affine algebraic variety, a Stein open set
$\Omega\subset Y$ is Runge if and only if the polynomial functions on~$Y$ are
dense in~$\oc(\Omega)$.

Our first result given below concerns approximations of holomorphic
maps by (complex) Nash algebraic maps. If $Y$, $Z$ are
quasi-projective (irreducible, reduced) algebraic varieties, a map
$f:\Omega\to Z$ defined on an open subset $\Omega\subset Y$ is said to be
{\sl Nash algebraic} if $f$ is holomorphic and the graph
$$\Gamma_f:=\big\{(y,f(y))\in\Omega\times Z: y\in\Omega\big\}$$
is contained in an algebraic subvariety $G$ of $Y\times Z$ of dimension
equal to~$\dim Y$. If~$f$ is Nash algebraic, then the image $f(\Omega)$
is contained in an algebraic subvariety $A$ of $Z$ with
$\dim A=\dim f(\Omega)\le\dim Y$. (Take $A=\pr_2(G)\subset Z$,
after eliminating any unnecessary components of~$G$.)

\than Theorem 1.1|Let
$\Omega$ be a Runge domain in an affine algebraic variety $S$, and let
$f:\Omega\to X$ be a holomorphic map into a quasi-projective algebraic
manifold $X$. Then for every relatively compact domain
$\Omega_0\compact\Omega$, there is a sequence of
Nash algebraic maps \hbox{$f_\nu:\Omega_0\to X$} such that
$f_\nu\to f$ uniformly on $\Omega_0$.

\noindent Moreover, if there is an
algebraic subvariety $A$ $($not necessarily reduced$)$ of $S$ and an algebraic
morphism $\alpha:A\to X$ such that $f_{|A\cap\Omega}= \alpha_{|A\cap\Omega}$,
then the $f_\nu$ can be chosen so that
$f_{\nu|A\cap\Omega_0}=f_{|A\cap\Omega_0}$. $($In particular, if we are given
a positive integer $k$ and a finite set of points $(t_j)$
in $\Omega_0$, then the $f_\nu$ can be taken to have the same $k$-jets
as $f$ at each of the points $t_j.)$ \fintha

An algebraic subvariety $A$ of $S$ is given by a coherent sheaf $\oc_A$ on
$S$ of the form $\oc_A=\oc_S/\ic_A$, where $\ic_A$ is an ideal
sheaf in $\oc_S$ generated by a (finite) set of polynomial functions on $S$.
(If $\ic_A$ equals the ideal sheaf of the algebraic subset
$\Supp\;\oc_A\subset S$, then $A$ is reduced and we can identify $A$ with
this subset.)  The restriction of a holomorphic map  $f:\Omega\to X$ to an
algebraic subvariety $A$ is given by
$$f_{|A\cap\Omega}=f\circ\iota_{A\cap\Omega}: A\cap\Omega\to X\;,$$ where
$\iota_{A\cap\Omega}:A\cap\Omega\to\Omega$ is the inclusion morphism.  We
note that the $k$-jet of a holomorphic map $f:\Omega\to X$ at a point
$a\in\Omega$ can be described as the restriction $f_{|\{a\}^k}:\{a\}^k\to X$
of $f$ to the nonreduced point $\{a\}^k$ with structure ring
$\oc_{\Omega,a}/{\cal M}^{k+1}_{\Omega,a}$ (where ${\cal M}_{\Omega,a}$
denotes the maximal ideal in $\oc_{\Omega,a}$).  The parenthetical statement
in Theorem~1.1 then follows by setting $A$ equal to the union of the
nonreduced points $\{t_j\}^k$.

If in Theorem~1.1 we are given an exhausting sequence $(\Omega_\nu)$ of
relatively compact open sets in $\Omega$, then we can construct by the Cantor
diagonal process a sequence of Nash algebraic maps
$f_\nu:\Omega_\nu\to X$ converging to $f$ uniformly on every compact subset
of $\Omega$.  We are unable to extend Theorem~1.1 to the case where $X$ is
singular.  Of course, for the case $\dim S=1$, Theorem~1.1 extends to
singular $X$, since then $f$ can be lifted to a
desingularization of $X$.

In the case where $X$ is an affine algebraic manifold, Theorem~1.1 is
easily proved as follows:  One approximates $f:\Omega\to X\subset\cb^m$ with
an algebraic map $g:\Omega_0\to\cb^m$ and then applies a Nash algebraic
retraction onto $X$ to obtain the desired approximations. (The
existence of Nash algebraic retractions is standard and is given in
Lemma~2.1.) In the case where $f$ is a map into a projective manifold
$X\subset\pb^{m-1}$, one can reduce to the case where $f$ lifts to a map $g$
into the cone $Y\subset\cb^m$ over $X$, but in general one cannot find a
global Nash algebraic retraction onto $Y$. Instead our proof proceeds in this
case by pulling back to $\Omega$ the normal bundle to $Y\ssm\{0\}$ (after
making some reductions) and then using this pull-back bundle to extend $g$ to
a map $G:\Omega\times\cb^p\to\cb^m$.  We then approximate $G$ by algebraic
maps $G_\nu$, which we compose (on the right) with Nash algebraic maps from
$\Omega_0\times\{0\}$ into $G_\nu^{-1}(Y)$ to obtain the approximations
$f_\nu$.

In the case where the map $f$ in Theorem~1.1 is an embedding of a smooth
domain $\Omega$ into a projective manifold $X$ with the property that $f^\star
L$ is trivial for some ample line bundle $L$ on $X$, we can approximate $f$ by
Nash algebraic embeddings $f_\nu$ with images contained in affine Zariski open
sets of the form $X\ssm D_\nu$, where the $D_\nu$ are divisors of powers of
$L$. This result (Theorem~4.1) is obtained by first modifying $f$ on a Runge
domain $\Omega_1\i\i\Omega$ in such a way as to create essential
singularities on the boundary of $f(\Omega_1)$ so that the closure
$f(\ov\Omega_1)$ becomes complete pluripolar (in the sense of
plurisubharmonic function theory). The existence of an ample divisor avoiding
$f(\Omega_0)$ is then obtained by means of an approximation theorem
(Proposition~4.8) based on H\"ormander's $L^2$ estimates for $\ov\partial$.

One of our main applications is the study of the Kobayashi pseudodistance
and the Kobayashi-Royden infinitesimal pseudometric on algebraic manifolds
([Ko67], [Ko70], [Ro71]). We use throughout the following notation:

\smallskip \noindent
$\Delta_R(a) = \{t\in \cb:|t-a|<R\}$,~ $\Delta_R = \Delta_R(0)$,~ $\Delta =
\Delta_1$.\hfil\break
For a complex manifold $X$:\hfil\break
$T_X={}$ the holomorphic tangent bundle of $X$.\hfil\break
$N_S=N_{S/X}={}$ the normal bundle $T_{X|S}/T_S$ of a submanifold
$S\subset X$.\hfil\break
$\kappa_X={}$ the Kobayashi-Royden infinitesimal Finsler pseudometric on
$T_X$,\hfil\break
\hglue 0.833cm
i.e., for $v\in T_X$,
$$\kappa_X(v)=\inf\Big\{\lambda>0\,:\,\exists
f:\Delta\to X~ \hbox{\rm holomorphic with}~
f_\star\Big(\lambda{\partial\over\partial\zeta}{}_{|0}\Big)=v\Big\}.$$
In addition, for a complex space $Z$, we let
$d_Z(a,b)$ denote the Kobayashi pseudodistance between points $a,b\in Z$.
In the case where $Z$ is smooth, then $d_Z(a,b)$ is the $\kappa_Z$-length of
the shortest curve from $a$ to $b$; see [NO90]. (A complex space $Z$ is said
to be {\sl hyperbolic} if $d_Z(a,b)>0$ whenever $a,b$ are distinct points of
$Z$.  A compact complex manifold is hyperbolic if and only if the
Kobayashi-Royden pseudometric is positive for all nonzero tangent vectors
[NO90].)

As an application of Theorem~1.1 (for the case where $\Omega$ is the unit
disk $\Delta$), we describe below how both the Kobayashi-Royden
pseudometric and the Kobayashi pseudodistance on projective algebraic
manifolds can be given in terms of algebraic curves.

\than Corollary 1.2|Let $M$ be an open subset of a projective
algebraic manifold $X$ and let $v\in T_M$.  Then $\kappa_M(v)$ is
the infimum of the set of $r\geq 0$ such that there exist a
closed algebraic curve $C\subset X$, a normalization $\gamma:C'\rightarrow
C\cap M$, and a vector $u \in T_{C'}$ with $\kappa_{C'}(u)=r$ and
$d\gamma(u)=v$. \fintha

\noindent {\sl Proof (assuming Theorem~1.1)}\pointir
Let $r_0$ be the infimum given above.  By the distance-decreasing
property of holomorphic maps, $\kappa_X (v)\leq r_0$.  We must verify
the reverse inequality.  Let $\varepsilon >0$ be arbitrary.  Choose a
holomorphic map $g:\Delta\rightarrow M$ and $\lambda\in
[0,\kappa_X(v)+\varepsilon]$ such that $dg(\lambda e_0)=v$,
where $e_0={\partial\over\partial\zeta}{}_{|0}$. By Theorem~1.1 with
$\Omega=\Delta$, there exists a Nash
algebraic map $f:\Delta _{1-\varepsilon}\rightarrow M$ such that
$df(e_0)=dg (e_0)$. Let $\wt f: \Delta\rightarrow M$ be given by
$\wt f(\zeta)=f((1-\varepsilon)\zeta)$ and let
$\wt\lambda=\lambda/(1-\varepsilon)$. Then $d\wt f(\wt\lambda e_0)=v$.
Since $\wt f$ is Nash algebraic, there exists an algebraic curve $C\subset X$
such that $\wt f(\Delta )\subset C$. Let $\gamma:C'\rightarrow C\cap M$ be the
normalization of $C\cap M$. We then have a commutative diagram
$$\eqalign{
\Delta&\buildo{{\displaystyle \wh f}}\over{\longrightarrow}C'\cr
\wt f&\searrow\quad\downarrow\gamma\cr
&\quad C\cap M\,.\cr}$$
Let $u = d\wh f(\wt \lambda e_0)$. We have $d\gamma(u)=v$. Thus
$$r_0\leq\kappa_{C'}(u) \leq\wt\lambda  = {\lambda \over 1-\varepsilon}\le
{\kappa_M(v)+\varepsilon\over 1-\varepsilon}\,.$$
Since $\varepsilon>0$ is
arbitrary, we conclude that $r_0\leq\kappa_M(v)$.\finpr\vskip7pt

If we let $M$ be Zariski-open in Corollary~1.2, we obtain the following
result:

\than Corollary 1.3|Let $Z$ be a quasi-projective algebraic manifold.
Then the Kobayashi-Royden pseudometric $\kappa_Z$ is given
by $$\kappa_Z(v)=\inf_C\kappa_{C'}(v),~~~~\forall v\in T_Z\,,$$
where $C$ runs over all $($possibly singular$)$ closed algebraic curves
in $Z$ tangent to~$v$, and $C'$ is the normalization of~$C$.
\fintha

In the case where $C'$ is compact and hyperbolic (i.e., has
genus $\geq 2$), the Kobayashi-Royden pseudometric $\kappa_{C'}$ just
coincides with the Poincar\'e metric induced by the universal covering
$\Delta\to C'\,$; thus in some sense the computation of
$\kappa_Z$ is reduced to a problem which is more algebraic in nature.
Of particular interest is the case where $Z$ is projective; then the
computation becomes one of determining the genus and Poincar\'e metric of
curves in each component of the Hilbert scheme. This observation strongly
suggests that it should be possible to characterize Kobayashi hyperbolicity
of projective manifolds in purely algebraic terms; this would follow for
instance from S.~Lang's conjecture [La86], [La87,~IV.5.7] that a projective
manifold is hyperbolic if and only if it has no rational curves and no
nontrivial images of abelian varieties; see also [De93] for further results
on this question.

We likewise can use algebraic curves to determine the Kobayashi
pseudodistance in quasi-projective algebraic varieties:

\than Corollary 1.4|Let $Z$ be a quasi-projective algebraic variety, and
let $a,b\in Z$.  Then the Kobayashi pseudodistance $d_Z(a,b)$ is given by
$$d_Z(a,b)=\inf_C d_C(a,b)\,,$$ where $C$ runs over all $($possibly
singular and reducible$)$ one-dimensional algebraic subvarieties of
$Z$ containing $a$ and $b$.\fintha

\rem Remark|If $C$ is a one-dimensional reducible complex space, then by
definition $d_C(a,b)$ is the infimum of the sums $$\sum_{j=1}^n
d_{C_j}(a_{j-1},a_j)$$ where $n$ is a positive integer,
$C_1,\ldots,C_n$ are (not necessarily distinct) irreducible components of
$C$, $a_0=a\in C_1$, $a_n=b\in C_n$ and  $a_j\in C_j\cap C_{j+1}$ for
$1\leq j\leq n-1$.

\noindent {\sl Proof of Corollary~1.4 (assuming Theorem~1.1)}\pointir Let
$\epsilon>0$ be arbitrary.  It suffices to find algebraic curves $C_j\subset
Z$ as in the above remark with $$\sum_{j=1}^n d_{C_j}(a_{j-1},a_j) < d_Z(a,b)
+\epsilon\;.$$  By the definition of the Kobayashi pseudodistance, there
exist points $a=a_0,a_1\ld a_n=b$ in $Z$ and holomorphic maps $f_j:\Delta\to
Z$ such that $f_j(0)=a_{j-1}, f_j(t_j)=a_j$ ($t_j\in\Delta$), for $1\leq
j\leq n$, and $\sum d_\Delta(0,t_j)<d_Z(a,b) +\epsilon$. Let $\rho<1$ be
arbitrary. By Theorem~1.1 (applied
to the desingularization of $Z$), we can find Nash algebraic maps
$g_j:\Delta_\rho\to Z$ with $g_j(0)=f_j(0)=a_{j-1}$ and
$g_j(t_j)=f_j(t_j)=a_j$.  Then for $1\leq j\leq n$, $g_j(\Delta)$ is
contained in an algebraic curve $C_j$.  Hence
$$\sum_{j=1}^nd_{C_j}(a_{j-1},a_j)\leq\sum_{j=1}^nd_{\Delta_\rho}(0,t_j)=
\sum_{j=1}^nd_\Delta\left(0,{t_j\over\rho}\right)\;.$$ Since $\rho<1$ is
arbitrary, our conclusion follows.\finpr\vskip7pt

Eisenman [Ei70] has also introduced an invariant pseudometric on the
decomposable $p$-vectors of an arbitrary complex space~$Y$. Let $B^p$ be
the unit ball of $\cb^p$. At each nonsingular point of~$Y$, the {\sl Eisenman
$p$-metric} $E^p_Y$ is defined by $$E^p_Y(v)=\inf\Big\{\lambda>0:\,\exists
f:B^p\to Y~ \hbox{\rm holomorphic
with}~f_\star\Big(\lambda{\partial\over\partial z_1}
\wedge\ldots\wedge{\partial\over\partial z_p} {}_{|0}\Big)=v\Big\}$$ for all
decomposable $p$-vectors $v\in\Lambda^pT_Y$. (The space $Y$ is said to be {\sl
strongly $p$-measure hyperbolic} if $E^p_Y(v)>0$ whenever $v\ne 0$.) For the
case $p=\dim Y$, $E^p_Y$ is called the {\sl Eisenman volume}. The following
generalization of Corollary~1.3 to the Eisenman $p$-metric is obtained by a
parallel argument.

\than Corollary 1.5|The Eisenman $p$-metric of a quasi-projective
algebraic manifold $Z$ can be computed in terms of the Eisenman volumes of
its $p$-dimensional algebraic subvarieties. More precisely
$$E^p_Z(v)=\inf_YE^p_Y(v)$$ for all decomposable $p$-vectors
$v\in\Lambda^pT_Z$, where $Y$ runs over all $($possibly singular$)$
$p$-dimensional closed algebraic subvarieties of~$Z$ tangent to~$v$.  $($A
decomposable $p$-vector $v$ at a singular point $y\in Y$ is said to be tangent
to $Y$ if the analytic germ of $Y$ at $y$ contains a smooth irreducible
component tangent to $v$; the above definition of $E^p_Y(v)$ carries over to
this case.$)$ \fintha

If $Y$ is compact and the canonical bundle $K_{\wt Y}$ of a
desingularization ${\wt Y}\to Y$ is ample, then a form of the maximum
principle tells us that the Eisenman volume $E^p_{\wt Y}$ is bounded from
below by the volume form of the K\"ahler-Einstein metric (with curvature
$-1$) on ${\wt Y}$ (see [Yau78], [NO90,~\S 2.5]). It would be very
interesting to know when these two intrinsic volumes are actually equal.  In
connection to this question, one would like to know whether the infimum
appearing in Corollary~1.5 can be restricted to subvarieties $Y$ with $\wt Y$
satisfying this property, in virtue of some density argument. This would make
the Eisenman $p$-metric theory completely parallel to the case of the
Kobayashi-Royden pseudometric $\kappa_Z=E^1_Z$.

The following result of Stout [St84] on the exhaustion of Stein manifolds by
Runge domains in affine algebraic manifolds allows us to extend Theorem~1.1
to maps from arbitrary Stein manifolds.

\than Theorem 1.6 {\rm(E.L.~Stout)}|Let $S$ be a Stein
manifold. For every Runge open set $\Omega\compact S$, there exists an
affine algebraic manifold $Y$ such that $\Omega$ is biholomorphic to a
Runge open set in~$Y$.
\fintha

Theorem~1.6 has been refined into a relative version by
Tancredi-Tognoli [TT89]. An immediate consequence of Theorem~1.1 and Stout's
Theorem~1.6 is the following approximation theorem for holomorphic maps from
Stein manifolds into projective algebraic manifolds:

\than Corollary 1.7|Let $\Omega$ be a relatively compact domain in a
Stein manifold $S$, and let $f:S\to X$ be a holomorphic map into a
quasi-projective algebraic manifold $X$.
Then there is a sequence of holomorphic maps \hbox{$f_\nu:\Omega\to X$}
such that $f_\nu\to f$ uniformly on $\Omega$ and the images $f_\nu(\Omega)$
are contained in algebraic subvarieties $A_\nu$ of $X$ with $\dim A_\nu=\dim
f_\nu(\Omega)$.  Moreover, if we are given
a positive integer $k$ and a finite set of points $(t_j)$
in $\Omega_0$, then the $f_\nu$ can be taken to have the same $k$-jets
as $f$ at each of the points $t_j$. \fintha

We must point out at this stage that Theorem~1.6 completely breaks
down for singular Stein spaces. In fact, there are examples of germs of
analytic sets with non-isolated singularities which are not
biholomorphic to germs of algebraic
sets; Whitney [Wh65] has even given an example which is not
$C^1$-diffeomorphic to an algebraic germ. The technique used by Stout
and Tancredi-Tognoli to prove Theorem~1.6 consists of a generalization
to the complex case of the methods introduced by Nash~[Nh52].
In the present work, we give a different proof of Theorem~1.6
based on some of the arguments used in the proof of the approximation
theorem~1.1.  Our proof starts with the observation that every Stein domain
$\Omega\compact S$ can be properly embedded in an affine algebraic manifold
in such a way that the normal bundle is trivial. This key step also yields the
following parallel result for holomorphic vector
bundles over Stein manifolds. (A similar result has also been given by
Tancredi and Tognoli [TT92, Theorem~4.1]; see Section~3.)

\than Theorem 1.8|Let $E$ be a holomorphic vector bundle on an
$n$-dimensional Stein manifold $S$. For every Runge open set $\Omega\compact
S$, there exist an $n$-dimensional projective algebraic manifold $Z$, an
algebraic vector bundle $\wt E\to Z$, an ample divisor $D$ in $Z$, and a
holomorphic injection $i:\Omega\hookrightarrow Z\ssm D$ with the
following properties:
\smallskip \item{\rm(i)} $i(\Omega)$ is a Runge domain in the affine
algebraic manifold $Z\ssm D\,;$ \vskip0pt
\item{\rm(ii)} $i^\star\wt E$ is isomorphic to~$E_{|\Omega}$. \fintha

Theorem~1.8 is proved by observing that there is a holomorphic map from
$S$ into a Grassmannian such that the given vector bundle $E$ on $S$
is isomorphic to the pull-back of the universal quotient vector bundle. The
desired result then follows from Theorem~1.6 applied to $\Omega\compact S$
and from the existence of Nash algebraic approximations of the map from
$\Omega$ to the Grassmannian (see Section~5).

A substantial part of the results of this work was obtained
during a stay of the third author at Institut Fourier, Universit\'e de
Grenoble~I, and he wishes to thank Institut Fourier for its hospitality.
\vskip12pt

\titre 2. Holomorphic and Nash algebraic retractions|

Let $S$ be an $n$-dimensional Stein manifold. The Bishop-Narasimhan
embedding theorem ([Bi61], [Na60]) implies that $S$ can be embedded as
a closed submanifold of a complex vector space~$\cb^N$ with e.g.
$N=2n+1\,$; by Eliashberg and Gromov [EG92], it is in fact enough to take
$N=[(3n+3)/2]$, and this value is optimal for even~$n$. Such optimal
values will not be needed here. We first recall a well-known and elementary
lemma about the existence of holomorphic retractions for Stein submanifolds.

\than Lemma 2.1|Let $M$ be a Stein $($not necessarily closed$)$ submanifold
of a complex manifold~$Y$ and let $n=\dim Y$.
\smallskip
\item{\rm(i)} There exist a neighborhood $U$ of the zero section in the
normal bundle $N_M$ of~$M$, a neighborhood $V$ of $M$ in $Y$ and a
biholomorphism $\psi:U\to V$, such that $\psi$ maps the zero section of~$N_M$
identically onto~$M$.
\smallskip
\item{\rm(ii)} Let $\pi:N_M\to M$ be the natural projection. Then
$\rho=\pi\circ\psi^{-1}:V\to M$ is a holomorphic retraction, i.e. a
holomorphic map such that $\rho_{|M}=\Id_M$. \smallskip
\noindent Moreover, if $Y$ is affine algebraic and $M$ is a closed algebraic
submanifold, then $\psi$ and
$\rho$ can be taken to be Nash algebraic.
\fintha

\proof We first outline the proof of the well-known case where $Y=\cb^n$.
(Details can be found in [GR65,~p.~257].) As $M$ is Stein, by Cartan's
Theorem~B the normal bundle sequence $$0\la T_M\la T_{\cb^n|M}\la N_M\la
0\eqno(2.1)$$ splits; i.e., there is a global morphism
$\sigma:N_M\to T_{\cb^n|M}$ which is
mapped to~$\Id_{N_M}$ by composition with the projection onto~$N_M$. Then
the map $\psi:N_M\to\cb^n$ given by
$$\psi(z,\zeta)=z+\sigma(z)\cdot\zeta\;,~~~~\zeta\in
N_{M,z}\;,~z\in M\;,\eqno(2.2)$$ coincides with $\Id_M$ on the zero section
(which we also denote by $M$) of $N_M$, and thus the derivative $d\psi$
satisfies the equality $(d\psi)_{|T_M}=\Id_{T_M}$. Furthermore, at each point
$z\in M$, $\psi$ has vertical derivative
$(d^v\psi)_z:=d\psi_{|N_{M,z}}=\sigma_z$, and hence $d\psi$ is invertible
at~$z$. By the implicit function theorem, $\psi$ defines a biholomorphism
from a neighborhood $U$ of the zero section of $N_M$ onto a neighborhood $V$
of~$M$ in~$\cb^n$. Then $\psi_{|U}:U\to V$ is the required biholomorphism, and
$\rho=\pi\circ(\psi_{|U})^{-1}$ is a retraction. When $M$ is affine
algebraic, $N_M$ can be realized as an affine algebraic manifold such that
$\pi:N_M\to M$ is an algebraic map. The splitting morphism $\sigma$ can also
be taken to be algebraic. Then $\Gamma_\psi$ is an open subset of an
$n$-dimensional subvariety $G\subset N_M\times Y$, and
$\Gamma_\rho=(\Id_Y\times\pi)(\Gamma_\psi^{-1})$ is contained in the
$n$-dimensional algebraic variety $(\Id_Y\times\pi)(G)$.

In the case of a general complex manifold~$Y$, it is enough to consider
the case where $Y$ is Stein, since by [Siu76] every Stein subvariety
of a complex space $Y$ has arbitrarily small Stein neighborhoods.
(See also [De90] for a simpler proof; this property is anyhow very easy to
check in the case of nonsingular subvarieties.) Thus we can
suppose $Y$ to be a closed $n$-dimensional submanifold of
some~$\cb^p$, e.g. with $p=2n+1$. Observe that the normal bundle
$N_{M/Y}$ of $M$ in $Y$ is a subbundle of $N_{M/\cb^p}$. By the first case
applied to the pairs $M\subset\cb^p$, $Y\subset\cb^p$, we get a splitting
$\sigma':N_M\to T_{\cb^p|M}$ of the normal bundle sequence for $N_{M/\cb^p}$,
a map $\psi':N_{M/\cb^p}\to \cb^p$ inducing a biholomorphism $\psi':U'\to
V'$, and a retraction $\rho'':V''\to Y$, where $U'$ is a neighborhood of the
zero section of~$N_{M/\cb^p}$, $V'$ a neighborhood of $M$ in $\cb^p$, and
$V''$ a neighborhood of~$Y$. We define $\psi$ to be the composition
$$\psi:N_{M/Y}\lhra N_{M/\cb^p}\buildo{\displaystyle \psi'}\over
{\relbar\joinrel\la}\cb^p\buildo{\displaystyle\rho''}\over
{\relbar\joinrel\la}Y.$$ Of course, $\psi$ is not defined on the whole
of~$N_{M/Y}$, but only on $N_{M/Y}\cap(\psi')^{-1}(V'')$. Clearly
$\psi_{|M}=\Id_M$, and the property $(d^v\psi)_{|M}={\sigma'}\!_{|N_{M/Y}}$
follows from the equalities $(d^v\psi')_{|M}=\sigma'$ and
$(d\rho'')_{|T_Y}=\Id_{T_Y}$. It follows as before that $\psi$ induces a
biholomorphism $U\to V$ from a neighborhood of the zero section of $N_{M/Y}$
onto a neighborhood $V$ of~$M$. If $Y$ and $M$ are algebraic, then we can
take $\psi'$ to be algebraic and $\Gamma_{\rho''}$ to be contained in an
algebraic $p$-dimensional subvariety $G''$ of $\cb^p\times Y$. Then
$\Gamma_\psi$ is contained in some $n$-dimensional component of the algebraic
set $$\big(N_{M/Y}\times Y\big)\cap(\psi'\times\Id_Y)^{-1}(G'').$$
Therefore $\psi$ is Nash algebraic, and $\rho=\pi\circ\psi^{-1}$ is
also Nash algebraic as before.\finpr

In fact, the retraction $\rho$ of Lemma~2.1 can be taken to be Nash
algebraic in slightly more general circumstances, as is provided by
Lemma~2.3 below. First we need the following elementary lemma due to
Tancredi and Tognoli:

\than Lemma 2.2 {\rm [TT89, Corollary 2]}|  Let $A$ be an algebraic
subvariety (not necessarily reduced) of an affine algebraic variety $S$,
and let $\Omega$ be a Runge domain in $S$.  Suppose that $h$ is a
holomorphic function on $\Omega$ and $\gamma$ is an algebraic function on
$A$ such that $h_{|A\cap\Omega}=\gamma_{|A\cap\Omega}$.  Then $h$ can be
approximated uniformly on each compact subset of $\Omega$ by algebraic
functions $h_\nu$ on $S$ with $h_{\nu|A}=\gamma$.\fintha\vskip7pt

\noindent Tancredi and Tognoli [TT89] consider only reduced subvarieties
$A$, but their proof is valid for nonreduced $A$.  We give the proof here
for the reader's convenience:  To verify Lemma~2.2, one first extends
$\gamma$ to an algebraic function $\wt\gamma$ on $S$.  By replacing $h$ with
$h-\wt\gamma$, we may assume that $\gamma=0$.  Choose global algebraic
functions $\tau_1\ld\tau_p$ generating the ideal sheaf $\ic_A$ at each
point of $S$, and consider the surjective sheaf homomorphism (of algebraic
or analytic sheaves) $\oc^p_S\buildo{\displaystyle\tau}\over{\la}\ic_A$
given by $\tau(f^1\ld f^p)=\sum\tau_jf^j$.  Now choose $f=(f^1\ld f^p)\in
H^0(\Omega,\oc^p)$ such that $\tau(f)=h$.  Since $\Omega$ is Runge in $S$,
we can approximate the $f^j$ by algebraic functions $f^j_\nu$ on $S$.  Then
the functions $h_\nu= \sum\tau_jf^j_\nu$ provide the desired algebraic
approximations of $f$.\vskip7pt

We now find it convenient to extend the usual definition of a retraction map
by saying that a map
$\rho:V\to S$ from a subset $V\i X$ to a subset $S\i X$ (where $X$ is any
space) is a retraction if $\rho(x)=x$ for all $x\in S\cap V$ (even if
$S\not\subset V$). The following lemma on approximating holomorphic
retractions by Nash algebraic retractions is a variation of a result of
Tancredi and Tognoli [TT90, Theorem~1.5]:

\than Lemma 2.3|Let $V$ be a Runge domain in $\cb^n$ and let $S$ be
an algebraic subset of $\cb^n$ such that the variety $M:=S\cap V$ is
smooth.  Suppose that $\rho:V\to S$ is
a holomorphic retraction and $V_0\compact V$ is a relatively compact domain
such that $\rho(V_0)\compact M$. Then there is a sequence of Nash
algebraic retractions $\rho_\nu:V_0\to M$ such that $\rho_\nu\to\rho$
uniformly on $V_0$.\fintha

\proof First we consider the special case where $\rho$ is the holomorphic
retraction constructed in Lemma~2.1. In this case, the conclusion follows from
the proof of Theorem~1.5 in [TT90], which is based on an ingenious
construction of Nash [Nh52].  We give a completely different short proof of
this case here: Let $V,S,M$ satisfy the hypotheses of Lemma~2.3, and suppose
that $\sigma:N_M\to T_{\cb^n|M},\;\psi:U\to V$, and
$\rho=\pi\circ\psi^{-1}:V\to M$ are given as in the proof of Lemma~2.1 (for
the case $Y=\cb^n$). We first construct a ``normal bundle sequence" for the
possibly singular subvariety $S$.  Let $f_1\ld f_r$ be polynomials generating
the ideal of $S$ at all points of $\cb^n$.  Let $$d\ic_S\subset\oc_S\otimes
T^\star_{\cb^n}=T^\star_{\cb^n|S}\simeq\oc^n_S$$ denote the coherent sheaf
on $S$ generated by the 1-forms $(1_{\oc_S}\otimes df_j)$. (This definition
is independent of the choice of generators $(f_j)$ of $\ic_S$.) We consider
the ``normal sheaf" $${\cal N}_S={\rm Image}\left(T_{\cb^n|S}\to{\cal
H}\!om(d\ic_S,\oc_S)\right)\;,$$ where we identify $T_{\cb^n|S}\simeq {\cal
H}\!om(T^\star_{\cb^n|S},\oc_S)$.  Thus we have an exact sequence of sheaves
of the form $$0\la{\cal K}\la
T_{\cb^n|S}\buildo{\displaystyle\tau}\over\la{\cal N}_S \la 0\;.\eqno(2.3)$$
Note that ${\cal N}_{S|M}=N_M$ and the restriction of (2.3) to the
submanifold $M$ is the exact sequence (2.1). Let $\fc={\cal H}\!om({\cal
N}_S,T_{\cb^n|S})$ and choose a surjective algebraic sheaf morphism
$\oc_S^p\buildo{\displaystyle\alpha}\over\to\fc$. Recall that $$\sigma\in
\Hom\;(N_M,T_{\cb^n|M})=H^0(M,\fc)\;.\eqno(2.4)$$ Since $M$ is Stein, we can
choose $h\in H^0(M,\oc_S^p)=\oc(M)^p$ such that $\alpha(h)=\sigma$.  Since
$V$ is Runge, we can find a sequence $h_\nu$ of algebraic sections in
$H^0(S,\oc_S^p)$ such that $h_\nu\to h$ uniformly on each compact subset of
$M$.  We let $$\sigma_\nu=\alpha(h_\nu)\in H^0(S,\fc)\;.\eqno(2.5)$$  Let
$\wt M\supset M$ denote the set of regular points of $S$. Then ${\cal
N}_{S|\wt M}=N_{\wt M}$ and the sections $\sigma_{\nu|\wt M}\in \Hom(N_{\wt
M},T_{\cb^n|\wt M})$ are algebraic. Recalling (2.2), we similarly define the
algebraic maps $\psi_\nu:N_{\wt M}\to \cb^n$ by
$$\psi_\nu(z,\zeta)=z+\sigma_\nu(z)\cdot\zeta\;,~~~~\zeta\in N_{\wt
M,z}\;,\;z\in\wt M\;.\eqno(2.6)$$  Choose open sets $V_1, V_2$ with
$V_0\compact V_1\compact V_2\compact V$. Then $\psi_\nu\to\psi$ uniformly
on the set $U'=\psi^{-1}(V_2)\compact U$,  and for $\nu\gg 0$,
$\psi_{\nu|U'}$ is an injective map with $\psi_\nu(U')\supset V_1$.  It
follows that the maps $\rho_\nu=\pi\circ\psi^{-1}_\nu:V_0\to M$ converge
uniformly on $V_0$ to $\rho$. (We remark that this proof gives algebraic
maps $\psi_\nu$, while the method of [TT90] provides only Nash algebraic
$\psi_\nu$.)

We now consider the general case. Choose a Runge domain $V_1$ such that
$V_0\compact V_1\compact V$ and $\rho(V_0)\compact V_1$, and let
$\varepsilon>0$ be arbitrary. It suffices to construct a Nash algebraic
retraction $\rho':V_0\to M$ with $\|\rho'-\rho\|<2\epsilon$ on $V_0$. By
Lemma~2.1 and the case already shown, there exist a neighborhood $W$ of
$M\cap V_1$ and a Nash algebraic retraction $R:W\to M\cap V_1$.  After
shrinking $W$ if necessary, we can assume that $\|R(z)-z\|<\varepsilon$
for all $z\in W$. By Lemma~2.2 (with $A,S$ replaced by $S,\cb^n$), we can
find a polynomial map $P=(P_1\ld P_n):\cb^n\to\cb^n$ such that
$\|P-\rho\|<\varepsilon$ on $V_0$, $P_{|S}=\Id_S$, and $P(V_0)\subset W$.
We consider the Nash algebraic retraction $$\rho'=R\circ P_{|V_0}: V_0\to
M\;.$$ Then $\|\rho'-\rho\|_{V_0}\leq \|R\circ P-P\|_{V_0}+
\|P-\rho\|_{V_0}<2\varepsilon$.\finpr

\than Lemma 2.4|Let $M$ be a Stein submanifold of a complex
manifold~$Y$,  and let $E\to Y$ be a holomorphic vector bundle. Fix
a holomorphic retraction \hbox{$\rho:V\to M$} of a neighborhood $V$
of $M$ onto~$M$. Then there is a neighborhood $V'\subset V$ of~$M$ such
that $E$ is isomorphic to $\rho^\star(E_{|M})$ on~$V'$. In particular,
if the restriction $E_{|M}$ is trivial, then $E$ is trivial on a
neighborhood of~$M$.\fintha

\proof By [Siu76] we can find a Stein neighborhood $V_0\subset V$. As
$V_0$ is Stein, the identity map from $E_{|M}$ to $(\rho^\star E)_{|M}$
can be extended to a section of $\Hom(E,\rho^\star E)$ over~$V_0$. By
continuity, this homomorphism is an isomorphism on a sufficiently small
neighborhood $V'\subset V_0$ of~$M$.\finpr\vskip7pt

Finally, we need an elementary lemma on the existence of generating
sections for holomorphic vector bundles on Stein spaces:

\than Lemma 2.5|Let $E$ be a holomorphic vector bundle of rank $r$
on an $n$-dimensional Stein space $S$. Suppose that $\ic$ is a coherent sheaf
of ideals in $\oc_S$, and let $A=\Supp\;\oc_S/\ic$. Then we can find $n+r$
holomorphic sections in $H^0(S,\ic\otimes E)$ generating $E$ at every point
of $S\ssm A$. \fintha

\rem Remark|Note that if $\ic=\oc_S$, then $A=\emptyset$. We also remark
that the precise bound $n+r$ for the number of generating sections is not
needed in this work.

\proof Cartan's Theorem B implies that there are global holomorphic sections
of $\ic\otimes E$ with prescribed values at any point or on any discrete
sequence of points of $S$.  If $g_1\ld g_k$ are holomorphic sections of $E$,
we write $$V_j(g_1\ld g_k)=\Bigl\{x\in S\ssm A:\dim \;\hbox{\rm
Span}\{g_1(x)\ld g_k(x)\}\leq j\Bigr\}\;,$$ for $0\leq j\leq r-1$. (Here we
use the convention $\dim\emptyset =-1$.) We shall construct $g_1\ld g_{n+r}$
in     $H^0(S,\ic\otimes E)$ such that $$\dim V_j(g_1\ld g_k)\leq
n-k+j~~~~\hbox{\rm for } 1\leq k\leq n+r\;,~~0\leq j\leq r-1\;.$$ (Note that
this inequality is vacuous for $j\geq k$.) In particular, setting $j=r-1$ and
$k=n+r$, we have $\dim V_{r-1}(g_1\ld g_{n+r})=-1$, and thus $\{g_1\ld
g_{n+r}\}$ spans $E$ at every point of $S\ssm A$.

We construct the $g_k$ inductively as follows: First, we choose a point in
each irreducible component of $S\ssm A$ and a section $g_1\in
H^0(S,\ic\otimes E)$ that is nonzero at each of these points; thus $\dim
V_0(g_1)\leq n-1$ as desired.  Now suppose that $g_1\ld g_k$ have been
given.  For $0\leq j\leq r-1$, we choose a discrete set $$A_j\subset
V_j(g_1\ld g_k)\ssm V_{j-1}(g_1\ld g_k)$$ such that $A_j$ contains one point
in each component of  $V_j(g_1\ld g_k)$ of the maximal dimension $n-k+j$.
(The sets $A_j$ may be finite, countably infinite, or empty; note that
$A_j=\emptyset$ for $j\geq k+1$.)  We now choose $g_{k+1}\in
H^0(S,\ic\otimes E)$ such that $$g_{k+1}(a)\not\in \hbox{Span}\{g_1(a)\ld
g_k(a)\}~~~\hbox{for } a\in A_1\cup\ldots\cup A_{r-1}\;,$$ and thus
$V_j(g_1\ld g_{k+1})\subset V_j(g_1\ld g_k)\ssm A_j$ for $0\leq j\leq r-1$.
It then follows from the choice of the sets $A_j$ that  $$\dim V_j(g_1\ld
g_{k+1})\leq n-k+j-1~~~~\hbox{\rm for } 0\leq j\leq r-1\;,$$ as desired.
\finpr

\than Corollary 2.6|Let $E$ be a holomorphic vector bundle over a finite
dimensional Stein space $S$. Then there is a holomorphic map $\Phi:S\to G$
into a Grassmannian such that $E$ is isomorphic to the pull-back
$\Phi^\star Q$ of the universal quotient vector bundle $Q$ on $G$.
\fintha

\proof By Lemma~2.5 (with $\ic=\oc_S$), we can find finitely many holomorphic
sections $g_1\ld g_m$ spanning $E$ at all points of $S$. We define the
holomorphic map $$\Phi:S\to G(m,r),~~~~x\mapsto
V_x=\left\{(\lambda_1\ld\lambda_m)\in\cb^m\,:\,\sum\lambda_jg_j(x)=0\right\}$$
into the Grassmannian $G(m,r)$ of subspaces of codimension $r$ in~$\cb^m$.
Then the generators $g_j$ define an isomorphism $\gamma:E\to\Phi^\star Q$
with $\gamma_x:E_x\to Q_{\Phi(x)}=\cb^m/V_x$ given by
$$\gamma_x(e)=\left\{(\mu_1\ld\mu_m)\in\cb^m:\sum\mu_jg_j(x)=e\right\}
\;,~~~e\in E_x\;,$$ for $x\in S$. (See, for example,
[GA74,~Ch.~V].) Hence $E\simeq\Phi^\star Q$.\finpr \vskip12pt

\bigtitle 3.  Nash algebraic approximation on Runge|
domains in affine algebraic varieties|

In this section, we prove the main approximation theorem~1.1 for
holomorphic maps on a Runge domain $\Omega$ in an affine algebraic
variety.  We also give an alternate proof and variation of a result of
Tancredi and Tognoli [TT93] on the exhaustion of holomorphic vector bundles
on $\Omega$ by algebraic vector bundles (Proposition~3.2).

We begin the proof of Theorem~1.1 by first making a reduction to the
case where $S=\cb^n$, $X$ is projective, and $f$ is an embedding.  To
accomplish this reduction, let $S$ be an algebraic subvariety of $\cb^n$ and
suppose that $\Omega\subset S$ and $f:\Omega\to X$ are given as in
Theorem~1.1.  By replacing $X$ with a smooth completion, we can assume that
$X$ is projective.  We consider the embedding $$f_1=(i_{\Omega},f):\Omega\to
X_1=\pb^n \times X,$$ where $i_{\Omega}$ denotes the inclusion map
$\Omega\hookrightarrow\cb^n \hookrightarrow \pb^n$. Approximations of $f_1$
composed with the projection onto~$X$ will then give approximations of $f$.
Since $\Omega$ is $\oc(S)$-convex, we can construct an analytic polyhedron
of the form $$\wt\Omega=\{z\in\cb^n :|h_j(z)|\leq 1\hbox{\ \ \rm for }1\leq
j\leq s\}\;,$$ where the $h_j$ are in $\oc(\cb^n)$, such that
$$\Omega_0\compact S\cap\wt\Omega\compact\Omega\;.$$ Then $\wt\Omega$ is a
Runge domain in $\cb^n$.  Consider the Runge domain
$$\wt\Omega_\varepsilon=\{z\in\wt\Omega:|P_j(z)|<\varepsilon\hbox{\ \ \rm
for }1\leq j\leq\ell\}\;,$$ where the $P_j$ are the defining polynomials for
$S\subset\cb^n$. In order to extend $f_1$ to $\wt\Omega_\varepsilon$, we
apply a result of Siu [Siu76] to find a Stein neighborhood $U\subset X'$ of
$f_1(\Omega)$. Then by using the holomorphic retraction of Lemma~2.1
(applied to $U$ embedded in an affine space), shrinking $\varepsilon$ if
necessary, we can extend $f_{1|S\cap\wt\Omega_\varepsilon}$ to a holomorphic
map $$\wt f_1:\wt\Omega_\varepsilon\to U\subset X_1\;.$$ By Nash
approximating the embedding $$f_2=(i_{\wt\Omega_\varepsilon},\wt
f_1):\wt\Omega_\varepsilon\to X_2=\pb^n\times X_1\;,$$ we also Nash
approximate $f$.

Thus by replacing $(\Omega,X,f)$ with $(\wt\Omega_\varepsilon,X_2,f_2)$,
we may assume that $\Omega$ is a Runge domain in $\cb^n$ and
$f:\Omega\to X$ is a holomorphic embedding. We
make one further reduction:  Choose a very ample line bundle $L$ on $X$.
By the proof of Lemma~2.5, we can find algebraic sections
$\sigma_1\ld\sigma_p$ of the algebraic line bundle $\alpha^\star L^{-1}\to
A$ without common zeroes. We then extend the sections $\sigma_{j|\Omega}$
to holomorphic sections $\wt\sigma_1\ld\wt\sigma_p$ of $f^\star L^{-1}$.
Again by Lemma~2.5, we choose holomorphic sections
$$\wt\sigma_{p+1}\ld\wt\sigma_{p+n+1}\in H^0(\Omega,\ic_A\otimes f^\star
L^{-1})$$ such that
$\wt\sigma_1\ld\wt\sigma_{p+n+1}$ generate $f^\star L^{-1}$ at all points of
$\Omega$.  By the proof of Corollary~2.6, the sections
$\wt\sigma_1\ld\wt\sigma_{p+n+1}$ define a holomorphic map $\Phi:\Omega\to
\pb^{p+n}$ and an isomorphism of line bundles $\gamma:f^\star
L^{-1}\to\Phi^\star\oc(1)$.  Similarly, the sections $\sigma_1\ld\sigma_p$
define an algebraic morphism $$\Phi_A:A\to\pb^{p-1}\subset\pb^{p+n}$$ and an
algebraic isomorphism $\gamma_A:\alpha^\star
L^{-1}\to\Phi_A^\star\oc(1)$ such that $\Phi_{|A\cap\Omega}=
\Phi_{A|A\cap\Omega}$ and $\gamma_{|A\cap\Omega}=\gamma_{A|A\cap\Omega}$.

We set $X'=\pb^{p+n}\times X$
and let $L'\to X'$ be the total tensor product $L'=\oc(1)\stimes L$. Then
$L'$ is very ample on $X'$ and is thus isomorphic to the hyperplane section
bundle of a projective embedding $X'\subset\pb^{m-1}$.  Hence $X'$ can be
identified with the projectivization of an affine algebraic cone
$Y\subset\cb^m$, and the line bundle $L'^{-1}=\oc_{X'}(-1)$ can be
identified with the blow-up $\wt Y$ of the cone $Y$ at its vertex. Let
$\tau:\wt Y\to Y$ denote the blow-down of the zero section $X'$ of $L'^{-1}$
(so that $\tau:\wt Y\ssm X'\approx Y\ssm\{0\}$).

We consider the maps $f'=(\Phi,f):\Omega\to X'~$,
${f'}\!_{|A\cap\Omega}={\alpha'}\!_{|A\cap\Omega}$. The pull-back bundle
$f^{\prime\star}L'$ is trivial; in fact, a trivialization of
$f^{\prime\star}L'$ is given by the nonvanishing global section $$\gamma\in
\Hom(f^\star L^{-1},\Phi^\star\oc(1))=H^0(\Omega,f^\star
L\otimes\Phi^\star\oc(1))=H^0(\Omega,f^{\prime\star}L')\;.$$  We can regard
the nonvanishing section $1/\gamma\in H^0(\Omega,f^{\prime\star}L'^{-1})$ as
a holomorphic mapping from $\Omega$ to $\wt Y\ssm X'$. Therefore
$f':\Omega\to X'$ lifts to the holomorphic embedding $$g=\tau\circ
(1/\gamma):\Omega\to Y\ssm\{0\}\;.$$  We likewise have an algebraic morphism
$$\beta=\tau\circ (1/\gamma_A):A\to Y\ssm\{0\}$$ such that $\beta$ is a lift
of $\alpha'$ and  $g_{|A\cap\Omega}=\beta_{|A\cap\Omega}$. It is sufficient
to Nash approximate $g$ by maps into $Y\ssm\{0\}$ agreeing with $\beta$ on
$A$. Hence Theorem~1.1 is reduced to the following:

\than Lemma 3.1|Let $\Omega$ be a Runge open set in $\cb^n$, let
$Y\subset\cb^m$ be an affine algebraic variety, and let $g:\Omega\to Y'$ be a
holomorphic embedding into the set $Y'$ of regular points of $Y$.  Suppose
that there is an algebraic subvariety $A$ (not necessarily reduced) of
$\cb^n$ and an algebraic morphism $\beta:A\to Y$ such that
$g_{|A\cap\Omega}=\beta_{|A\cap\Omega}$.  Then $g$ can be approximated
uniformly on each relatively compact domain \hbox{$\Omega_0\compact$}
$\Omega$ by Nash algebraic embeddings $g_\nu:\Omega_0\to Y'$ such that
$g_{\nu|A\cap\Omega_0}=  g_{|A\cap\Omega_0}$.\fintha

\proof We write $$Z=g(\Omega)\subset Y'\subset Y\subset \cb^m\;,$$ and we
let $N_Z$ denote the normal bundle of $Z$ in $\cb^m$. By the argument at the
beginning of the proof of Lemma~2.1, since $\Omega$ is Stein, the bundle
sequence (on $\Omega$) $$0\la g^\star T_{Z}\la g^\star T_{\cb^m}\la g^\star
N_{Z}\la 0\eqno(3.1)$$ splits, so that we obtain a holomorphic subbundle $N$
of $g^\star T_{\cb^m}$ with the property that $$g^\star T_{\cb^m}=N\oplus
g^\star T_{Z}\;.\eqno(3.2)$$  Identifying $g^\star T_{\cb^m}$ with
$\Omega\times\cb^m$, we define a map $\varphi :N\to\cb^m$ by
$$\varphi(z,w)=g(z)+w\hbox{\ \ \rm for } (z,w)\in N\subset g^\star
T_{\cb^m}=\Omega\times\cb^m\;.\eqno(3.3)$$ As in
the proof of Lemma~2.1, $\varphi$ defines
a biholomorphism from a neighborhood $\Omega'$ of the zero section
$\Omega\times\{0\}$ of $N$ onto a neighborhood $W$ of $Z$ in $\cb^m$ with the
property that $\varphi(z,0)=g(z)$, for $z\in\Omega$.  We can assume that $W$
does not contain any singular points of $Y$.

By Lemma~2.5 we
can choose sections $s_1,\ldots ,s_p$ of $N$ that span at every point of
$\Omega$.  Let
$\lambda:\Omega\times\cb^p \to N$ be the surjective vector-bundle homomorphism
given by $$\lambda(z,v_1,\ldots ,v_p)=\sum_{j=1}^p v_js_j(z)\;,\eqno(3.4)$$
for $z\in\Omega$, $v_1,\ldots,v_p\in\cb$, and consider the map
$$G=\varphi\circ\lambda:\Omega\times\cb^p\to\cb^m\;.$$
Let $\wt\Omega =\lambda^{-1}(\Omega')\supset \Omega\times\{0\}$.
Since $\varphi_{|\Omega'}$ is biholomorphic and
$\lambda$ is a submersion, it follows that $G_{|\wt\Omega}:\wt\Omega\to W$ is
also a submersion. Let $M=G^{-1}(Y)\cap \wt\Omega$.  Since $G_{|\wt\Omega}$
is a submersion and $Y\cap W=Y'\cap W$, it follows that
$M$ is a closed submanifold of $\wt\Omega$. We note that $G(z,0)=g(z)$ for
$z\in\Omega$, and thus $M \supset\Omega\times\{0\}$.

Choose Runge domains $\Omega_1,\Omega_2$ in $\cb^n$ with
$\Omega_0\compact\Omega_1\compact\Omega_2\compact\Omega$ and consider the
domains $\wt\Omega_1\compact\wt\Omega_2$ in $\cb^{n+p}$ given by
$\wt\Omega_i=\Omega_i\times\Delta_{i\varepsilon}^p$, where
$\varepsilon >0$ is chosen so that $\wt\Omega_2\compact\wt\Omega$. By the
proof of Lemma~2.1, there is a neighborhood $V\subset\wt\Omega$ of
$M$ together with a biholomorphism $\psi:U\to V$, where $U\i
N_M$ is a neighborhood of the zero section
$M$, such that $\psi_{|M}=\Id_M$. Let
$$V_i=\{z\in\wt\Omega_i:|h_j(z)|< i\delta\hbox{\ \ \rm for }1\leq
j\leq\ell\}\;,~~~i=1,2\;,$$ where we choose $h_1\ld h_\ell\in\oc(\wt\Omega)$
defining $M$, and we choose $\delta>0$ such that
$V_2\compact V$. We then have the holomorphic retraction
$\rho=\pi\circ\psi^{-1}:V\to M$, where $\pi:N_M\to M$ is the projection. By
shrinking $\delta$ if necessary, we can assume that $\rho(V_1)\i M\cap
V_2=M\cap\wt\Omega_2$.

We identify $A\subset\cb^n$ with $A\times\{0\}\subset\cb^{n+p}$
so that we can regard $A$ as an algebraic subvariety (not necessarily
reduced) of $\cb^{n+p}$; then $G_{|A}=\beta$.  Since $\wt\Omega$ is Runge, by
Lemma~2.3 we can choose a sequence of polynomial maps
$G_\nu:\cb^{n+p}\to\cb^m$ such that $G_{\nu|A}=\beta$ and $G_\nu\to G$
uniformly on a Runge domain $D\i\wt\Omega$ containing $\overline{V_2}$.  We
write $M_\nu=G_\nu^{-1}(Y)\cap D'$, where $D'$ is a Runge domain with
$V_2\compact D'\compact D$.  Then for $\nu\gg 0$, $M_\nu$ is a submanifold
of $D'$. Moreover, since $\pi_{|M}$ is the identity, for large $\nu$ the
restriction of $\pi$ to $\psi^{-1}(M_\nu)\cap\pi^{-1}(M\cap V_2)$ is a
biholomorphic map onto $M\cap V_2$.  The inverse of this latter map is a
section $t_\nu\in H^0(M\cap V_2,N_M)$, and $t_\nu\to 0$ uniformly on $M\cap
V_1$. We now consider the holomorphic retractions $$\rho_\nu=\psi\circ
t_\nu\circ\rho_{|V_1}:V_1\to M_\nu\;.\eqno(3.5)$$ Then $\rho_\nu\to\rho$
uniformly on  $V_1$.  Choose a domain $V_0$ such that
$\overline{\Omega_0}\times\{0\}\i V_0\compact V_1$ and $\rho(V_0)\compact
V_1$.  Since $V_1$ is Runge in $\cb^{n+p}$ and $\rho_\nu(V_0)\compact
V_1$ for $\nu\gg 0$, by Lemma~2.3 we can find Nash
algebraic retractions $\rho'_\nu:V_0\to M_\nu$ sufficiently close to
$\rho_\nu$ such that $\rho'_\nu\to\rho$ uniformly on  $V_0$.  Then the Nash
algebraic approximations $g_\nu:\Omega_0\to Y'$  of $g$ can be given by
$$g_\nu(z)=G_\nu\circ\rho'_\nu(z,0)\;,\;z\in\Omega_0\;.\eqno(3.6)$$

It remains to verify that
$$g_{\nu|A\cap\Omega_0}=\beta_{|A\cap\Omega_0}\;.\eqno(3.7)$$  Since
$\rho'_\nu$ is a retraction to $M_\nu$, we have
$$G_\nu\circ\rho'_{\nu|M_\nu}=G_\nu\;.\eqno(3.8)$$  It suffices to show that
$\ic_{M_\nu}\i\ic_{A|D'}$ (where $\ic_A$ now denotes the ideal sheaf on
$\cb^{n+p}$ with $\oc_A=\oc_{\cb^{n+p}}/\ic_A$), since this would give
us an inclusion morphism $A\cap D'\hookrightarrow M_\nu$, and then (3.7) would
follow by restricting (3.8) to $A\cap D'$ and recalling that
$G_{\nu|A}=\beta$. Let $h\in \ic_{Y,G_\nu(a)}$ be arbitrary, where $a\in
A\cap D'$. Since $G_{\nu|A}=G_{|A}$, it follows that $$h\circ G_\nu-h\circ
G\in\ic_{A,a}\;.\eqno(3.9)$$  Since $\Omega\times\{0\}\i M$, we have $$h\circ
g\in\ic_{M,a}\i\ic_{\Omega\times\{0\},a}\i\ic_{A,a}\;,$$ and thus it follows
from (3.9) that $h\circ G_\nu\in\ic_{A,a}$. Since
$G_\nu^\star\ic_{Y,G_\nu(a)}$ generates $\ic_{M_\nu,a}$, it follows that
$\ic_{M_\nu,a}\i\ic_{A,a}$.\finpr\vskip7pt

Tancredi and Tognoli ([TT93],~Theorem~4.1) showed that if $E$ is a
holomorphic vector bundle on a Runge domain $\Omega$ in an affine algebraic
variety $S$, then for every  domain $\Omega_0\i\i\Omega$, $E_{|\Omega_0}$
is isomorphic to a ``Nash algebraic vector bundle" $E'\to\Omega_0$.
Tancredi and Tognoli demonstrate this result by applying a result of Nash
[Nh52] on Nash-algebraically approximating analytic maps into the
algebraic manifold of $n\times n$ matrices of rank exactly $r$.
(Alternately, one can approximate analytic maps into a Grassmannian.)
Theorem~1.1 is a generalization of (the complex form of) this
approximation result of Nash. In fact, Theorem~1.1 can be used to obtain
the following form of the Tancredi-Tognoli theorem with a more explicit
description of the equivalent Nash algebraic vector bundle:

\than Proposition 3.2|Let $\Omega$ be a Runge domain in an $n$-dimensional
affine algebraic variety $S$ and let $E\to\Omega$ be a holomorphic vector
bundle. Then for every relatively compact domain $\Omega_0\compact\Omega$
there exist:
\smallskip
\item{\rm(i)} an $n$-dimensional projective algebraic variety $Z$,
\item{\rm(ii)}a Nash algebraic injection $i:\Omega_0\hookrightarrow Z$,
\item{\rm(iii)} an algebraic vector bundle $\wt E\to Z$ such that $i^\star\wt
E$ is isomorphic to $E_{|\Omega_0}$,
\item{\rm(iv)}an ample line bundle $L\to Z$ with trivial restriction $i^\star
L$. \smallskip\noindent Moreover, if $S$ is smooth, then $Z$ can be taken to
be smooth. \fintha

\proof By Corollary~2.6, there is a holomorphic map $\Phi:\Omega\to G$
into a Grassmannian, such that $E\simeq\Phi^\star Q$ where $Q$ is the
universal quotient vector bundle. We construct a holomorphic map of the form
$$f=(\Phi,\Phi'):\Omega\to X= G\times\pb^N,$$ where
$\Phi':S\to\pb^N$ is chosen in such a way that $f$ is an embedding and $X$
has an ample line bundle $L$ with trivial pull-back $f^\star L$ (using the
same argument as in the reduction of Theorem~1.1 to Lemma~3.1).  Then
$E\simeq f^\star Q_X$, where $Q_X$ is the pull-back of $Q$ to~$X$.

Choose a Runge domain $\Omega_1$ with
$\Omega_0\compact\Omega_1\compact\Omega$. By Theorem~1.1, $f_{|\Omega_1}$ is a
uniform limit of Nash algebraic embeddings $f_\nu:\Omega_1\to X$. Hence the
images $f_\nu(\Omega_1)$ are contained in algebraic subvarieties
$A_\nu\subset X$ of dimension $n=\dim S$.  Let $f'=f_{\mu|\Omega_0}$, where
$\mu$ is chosen to be sufficiently large so that $f'$ is homotopic
to~$f_{|\Omega_0}$.  Then the holomorphic vector bundle $E':=f'^\star Q_X$ is
topologically isomorphic to~$E_{|\Omega_0}$. By~a theorem of Grauert [Gr58],
$E'\simeq E_{|\Omega_0}$ as holomorphic vector bundles. If $S$ is
singular, we take $Z=A_\mu, i=f'$, and $\wt E=E'$.

If $S$ is smooth, we must modify $A_\mu$ to obtain an appropriate smooth
variety $Z$. To accomplish this, let $\sigma:A''_\mu\to A_\mu$ be the
normalization of~$A_\mu$ and let $f'':\Omega_0\to A''_\mu$ be the map such
that $f'=\sigma\circ f''\,$; since $f'$ is an embedding into~$X$, we see that
$f''(\Omega_0)$ is contained in the set of regular points of~$A''_\mu$. (The
reason for the introduction of the normalization is that $f'(\Omega_0)$ can
contain multiple points of~$A_\mu$.) By Hironaka's resolution of
singularities [Hi64], there exists a desingularization $\pi:Z\to A''_\mu$
with center contained in the singular locus of~$A''_\mu$. Let $i:\Omega_0
\hookrightarrow Z$ denote the embedding given by $f''=\pi\circ i$. We note
that the exceptional divisor of $Z$ does not meet $i(\Omega_0)$. The
algebraic vector bundle $\wt E:=(\sigma\circ\pi)^\star Q_X$ then satisfies
$$i^\star\wt E=(\sigma\circ\pi\circ i)^\star Q_X=f^{\prime\star}Q_X =E'\simeq
E_{|\Omega_0}.$$

Finally, we note that the line bundle $L'':=\sigma^\star L_{|A_\mu}$ is ample
on $A''_\mu$ and thus there is an embedding $A''_\mu\subset\pb^m$ such that
$\oc(1)_{|A''_\mu}=L^b$ for some $b\in\nb$ (where $\oc(1)$ is the hyperplane
section bundle in $\pb^m$).  We can suppose that $\pi$ is an embedded
resolution of singularities with respect to this embedding; i.e., there is a
modification $\pi:\wt X\to\pb^m$ such that $\wt X$ is smooth and $Z\subset\wt
X$ is the strict transform of $A''_\mu$. Then there is an exceptional divisor
$D\geq 0$ of the modification $\pi:\wt X\to\pb^m$ such that
$\wt L:=\pi^\star\oc(a)\otimes\oc(-D)$ is ample on $\wt X$ for $a\gg 0$.
(We can take $D=\pi^\star[\pi(H)]-[H]$, for an ample hypersurface
$H\subset\wt X$ that does not contain any component of the exceptional divisor
of $\wt X$.) We can assume that $\mu$ has been chosen sufficiently large so
that $f'^\star L\simeq (f^\star L)_{|\Omega_0}$ is trivial.  By construction
$i^\star\oc_Z(D)$ is trivial, and hence so is $i^\star(\wt L_{|Z})$.\finpr

\rem Remark 3.3|Here, it is not necessary to
use Grauert's theorem in its full strength. The special case we need
can be dealt with by means of the results obtained in Section~2. In
fact, for $\mu$ large, there is a one-parameter family of holomorphic
maps $F:\Omega_1\times U\to X$, where $U$ is a neighborhood of the
interval $[0,1]$ in~$\cb$, such that $F(z,0)=f(z)$ and $F(z,1)=f_\mu(z)$
on~$\Omega_1\,$; this can be seen by taking a Stein neighborhood $V$ of
$f(\Omega)$ in~$X$ and a biholomorphism $\varphi=\psi^{-1}$ from a
neighborhood of the diagonal in $V\times V$ onto a neighborhood of the
zero section in the normal bundle (using Lemma~2.1$\,$(i)); then
$F(z,w)=\psi\big(w\,\varphi(f(z),f_\mu(z))\big)$ is the required family.
Now, if $Y$ is a Stein manifold and $U$ is connected, Lemma~2.4 and
an easy connectedness argument imply that all slices over $Y$ of a
holomorphic vector bundle $B\to Y\times U$ are isomorphic, at least
after we take restriction to a relatively compact domain in~$Y$. Hence
$$E'=f_\mu^\star Q_{X|\Omega_0}=(F^\star Q_X)_{|\Omega_0\times\{1\}}
\simeq(F^\star Q_X)_{|\Omega_0\times\{0\}}=f^\star Q_{X|\Omega_0}
\simeq E_{|\Omega_0}.$$

\rem Remark 3.4|In general, the embedding $i:\Omega_0\hookrightarrow Z$
in Proposition~3.2 will not extend to a (univalent) rational map
$S\merto Z$, but will instead extend as a multi-valued branched map.
For example, let $S=X\ssm H$, where $H$ is a smooth hyperplane section
of a projective algebraic manifold $X$ with some nonzero Hodge numbers
off the diagonal.  By Cornalba-Griffiths [CG75,~\S~21] (see also
Appendix~2 in [CG75]) there is a holomorphic vector bundle $E'\to S$
such that the Chern character $\ch(E')\in H^{\rm even}(S,\qb)$ is not
the restriction of an element of $H^{\rm even}(X,\qb)$ given by
(rational) algebraic cycles.  Choose $\Omega_0\compact\Omega\subset S$
such that the inclusion $\Omega_0\hookrightarrow S$ induces an
isomorphism $H^{\rm even}(S,\qb)\approx H^{\rm even}(\Omega_0,\qb)$,
and let $E=E'_{|\Omega}$. Now let $i:\Omega_0\hookrightarrow Z$ and
$\wt E\to Z$ be as in the conclusion of Proposition~3.2.  We assert
that $i$ cannot extend to a rational map on $S$. Suppose on the
contrary that $i$ has a (not necessarily regular) rational extension
$i':X\merto Z$.  By considering the graph $\Gamma_{i'}\subset X\times
Z$ and the projections $\pr_1:\Gamma_{i'}\to X$, $\pr_2:\Gamma_{i'}\to Z$,
we obtain a coherent algebraic sheaf $\fc=\pr_{1\star}
\oc(\hbox\pr_2^\star \wt E)$ on $X$ with restriction $\fc_{|\Omega_0}
=i^\star\wt E\simeq E_{|\Omega_0}$, and therefore $\ch({\cal
F})_{|\Omega_0}=\ch(E')_{|\Omega_0}$. Now $\ch(\fc)$ is given by
algebraic cycles in $X$, but by the above, $\ch(\fc)_{|S}=\ch(E')$,
which contradicts our choice of $E'$. In particular,
$E_{|\Omega_0}$ is not isomorphic to an algebraic vector bundle on~$S$.

\rem Problem 3.5|The method of proof used here leads to the following
natural question. Suppose that $f:\Omega\to X$ is given as in Theorem~1.1, and
assume in addition that $S$ is nonsingular and $\dim X\ge 2\dim S+1$.
Is it true that $f$ can be approximated by Nash algebraic embeddings
$f_\nu:\Omega_\nu\to X$ such that $f_\nu(\Omega_\nu)$ is contained in
a {\sl nonsingular} algebraic subvariety $A_\nu$ of $X$ with
$\dim A_\nu=\dim f_\nu(\Omega_\nu)=\dim S\,$? If~$\dim X=2\dim S$, does
the same conclusion hold with $f_\nu$ being an immersion and $A_\nu$ an
immersed nonsingular variety with simple double points$\,$?

Our expectation is that the answer should be positive in both cases, because
there is a priori enough space in~$X$ to get the smoothness of
$f_\nu(\Omega_\nu)$ by a generic choice of~$f_\nu$, by a Whitney-type
argument. The difficulty is to control the singularities introduced by the
Nash algebraic retractions. If we knew how to do this, the technical point
of using a desingularization of $A_\mu$ in the proof of Proposition~3.2 could
be avoided.\vskip12pt

\titre 4. Nash algebraic approximations omitting ample divisors|

In this section we use $L^2$ approximation techniques to give the following
conditions for which the Nash algebraic approximations in
Theorem~1.1 can be taken to omit ample divisors:

\than Theorem 4.1|Let
$\Omega$ be a Runge domain in an affine algebraic manifold $S$, and let
$f:\Omega\to X$ be a holomorphic embedding into a projective algebraic
manifold $X$ $($with $\dim S<\dim X)$. Suppose that there exists an ample
line bundle $L$ on~$X$ with trivial restriction $f^\star L$ to~$\Omega$. Then
for every relatively compact domain $\Omega_0\compact\Omega$, one can find
sections $h_\nu\in H^0(X,L^{\otimes\mu(\nu)})$ and Nash algebraic embeddings
$f_\nu:\Omega_0 \to X\ssm h_\nu^{-1}(0)$ converging uniformly to
$f_{|\Omega_0}$. Moreover, if we are given a positive integer $m$ and a finite
set of points $(t_j)$ in $\Omega_0$, then the $f_\nu$ can be taken to have
the same $m$-jets as $f$ at each of the points $t_j$.\fintha

\than Corollary 4.2|Let $\Omega_0\compact\Omega\subset S$ be as in
Theorem~$4.1$. Suppose that $H^2(\Omega; \zb)=0$. Then for every
holomorphic embedding $f:\Omega\to X$ into a projective algebraic manifold
$X$ $($with $\dim S<\dim X)$, one can find affine Zariski open sets
$Y_\nu\subset X$ and Nash algebraic embeddings $f_\nu :\Omega_0\to Y_\nu$
converging uniformly to $f_{|\Omega_0}$.  Moreover, if we are given a
positive integer $m$ and a finite set of points $(t_j)$ in $\Omega_0$, then
the $f_\nu$ can be taken to have the same $m$-jets as $f$ at each of the
points $t_j$.\fintha

The main difficulty in proving Theorem~4.1 is that the submanifold $f(\Omega)$
might not be contained in an affine open set. To overcome this difficulty, we
first shrink $f(\Omega)$ a little bit and apply a small perturbation to make
$f(\Omega)$ smooth up to the boundary with essential singularities at each
boundary point. Then $f(\Omega)$ becomes in some sense very far from being
algebraic, but with an additional assumption on $f(\Omega)$ it is possible
to show that $f(\Omega)$ is actually contained in an affine Zariski open
subset (Corollary~4.9). The technical tools needed are H\"ormander's $L^2$
estimates for $\ov\partial$ and the following notion of complete pluripolar
set.

\than Definition 4.3|Let $Y$ be a complex manifold.
\smallskip
\item{\rm(i)} A function $u:Y\to[-\infty,+\infty[$ is said to be
quasi-plurisubharmonic $($quasi-psh for short$)$ if $u$ is locally equal to
the sum of a smooth function and of a plurisubharmonic (psh) function,
or equivalently, if $\sqrt{-1}\partial\op u$ is bounded below by a continuous
real $(1,1)$-form. \smallskip \item{\rm(ii)} A closed set $P$ is said to be
complete pluripolar in~$Y$ if there exists a quasi-psh function $u$ on $Y$
such that $P=u^{-1}(-\infty)$.
\smallskip\fintha

\rem Remark 4.4|Note that the sum and the maximum of two quasi-psh
functions is quasi-psh.  (However, the decreasing limit of a sequence
of quasi-psh functions may not be quasi-psh, since any continuous
function is the decreasing limit of smooth functions and there are many
continuous non-quasi-psh functions, e.g., $-\log^+|z|$.)  The usual
definition of pluripolar sets deals with psh functions rather than
quasi-psh functions; our choice is motivated by the fact that we want
to work on compact manifolds, and of course, there are no nonconstant
global psh functions in that case.  It is elementary to verify that if
$u$ is a quasi-psh function on a Stein manifold $S$, then $u+\varphi$
is psh for some smooth, rapidly growing function $\varphi$ on $S$, and
hence our definition of complete pluripolar sets coincides with the
usual definition in the case of Stein manifolds. Finally, we remark
that the word ``complete'' refers to the fact that $P$ must be the full
polar set of~$u$, and not only part of it.

\than Lemma 4.5|Let $P$ be a closed subset in a complex manifold~$Y$.
\smallskip
\item{\rm(i)} If $P$ is complete pluripolar in~$Y$,
there is a quasi-psh function $u$ on $Y$ such that $P=u^{-1}(-\infty)$
and $u$ is smooth on $Y\ssm P$.
\smallskip
\item{\rm(ii)} $P$ is complete pluripolar in $Y$ if and
only if there is an open covering $(\Omega_j)$ of $Y$ such that
$P\cap\Omega_j$ is complete pluripolar in~$\Omega_j$.
\smallskip
\fintha

\proof We first prove result (i) locally, say on a ball $\Omega=B(0,r_0)
\subset\cb^n$, essentially by repeating the arguments given in Sibony
[Sib85]. Let $P$ be a closed set in~$\Omega$, and let
$v$ be a psh function on $\Omega$ such that $P=v^{-1}(-\infty)$.
By shrinking $\Omega$ and subtracting a constant from $v$, we may assume
$v\leq 0$. Select a convex increasing function $\chi:[0,1]\to\rb$ such that
$\chi(t)=0$ on $[0,1/2]$ and $\chi(1)=1$. We set
$$w_k=\chi\big(\exp(v/k)\big).$$
Then $0\leq w_k\le1$, $w_k$ is plurisubharmonic on $\Omega$ and
$w_k=0$ in a neighborhood of~$P$. Fix a family $(\rho_\varepsilon)$
of smoothing kernels and set $\Omega'=B(0,r_1)\compact\Omega$.
For each~$k$, there exists $\varepsilon_k>0$ small such that
$w_k\star\rho_{\varepsilon_k}$ is well defined on $\Omega'$ and vanishes on
a neighborhood of $P\cap\Omega'$. Then
$$h_k=\max\big\{w_1\star\rho_{\varepsilon_1},\ldots,w_k\star
\rho_{\varepsilon_k}\big\}$$
is an increasing sequence of continuous psh functions such that
$0\leq h_k\leq 1$ on $\Omega'$ and $u_k=0$ on a neighborhood of
$P\cap\Omega'$.  The inequalities
$$h_k\geq w_k\star\rho_{\varepsilon_k}\geq w_k=\chi\big(\exp(v/k)\big)$$
imply that $\lim h_k=1$ on $\Omega'\ssm P$. By Dini's lemma, $(h_k)$
converges uniformly to~$1$ on every compact subset of~$\Omega'\ssm P$.
Thus, for a suitable subsequence $(k_\nu)$, the series
$$h(z):=|z|^2+\sum_{\nu=0}^{+\infty}(h_{k_\nu}(z)-1)\eqno(4.1)$$
converges uniformly on every compact subset of $\Omega'\ssm P$ and defines
a strictly plurisubharmonic function on $\Omega'$ which is continuous on
$\Omega'\ssm P$ and such that $h=-\infty$ on $P\cap\Omega'$. Richberg's
regularization theorem [Ri68] implies that there is a psh function $u$
on $\Omega'$ such that $u$ is smooth on $\Omega'\ssm P$ and $h\leq u\leq h+1$.
Then $u^{-1}(-\infty)=P\cap\Omega'$ and Property (i) is proved
on~$\Omega'$.

(ii) The ``only if'' part is clear, so we just consider the ``if'' part.
By means of (i) and by taking a refinement of the covering,
we may assume that we have open sets
$\wt\Omega_j\supset\!\supset\Omega_j\supset\!\supset\Omega'_j$ where
$(\Omega'_j)$ is a locally finite open
covering of~$Y$ and for each $j$ there is a psh function $u_j$
on $\wt\Omega_j$ such that
$P\cap\wt\Omega_j=u_j^{-1}(-\infty)$ with $u_j$ smooth on~$\wt\Omega_j\ssm P$.
We define inductively a strictly decreasing sequence
$(t_\nu)$ of real numbers by setting $t_0=0$ and
$$t_{\nu+1}\leq\inf\big\{\bigcup_{j,k\leq\nu}
u_j\big(\,\ov\Omega_j\cap\ov\Omega_k\ssm
u_k^{-1}([-\infty,t_\nu[)\big)\big\}\;.\eqno(4.2)$$
The infimum is finite since $u_j$ and $u_k$ have the same poles. Now,
if we take in addition $t_{\nu+1}-t_\nu<t_\nu-t_{\nu-1}$, there
exists a smooth convex function $\chi:\rb\to\rb$ such that
$\chi(t_\nu)=-\nu$.  Choose functions $\varphi_j\in C^\infty(\Omega_j)$
such that $\varphi_j \equiv 0$ on $\Omega'_j$,\break
$\varphi_j \leq 0$ on $\Omega_j$, and $\varphi_j(z)\to -\infty$ as $z
\to \partial\Omega_j$.  We set $$v_j=\chi\circ u_j+\varphi_j\;.$$
Then $v_j$ is quasi-psh on $\Omega_j$.  Let
$$v(z)=\max_{\Omega_j\ni z}~~v_j(z),~~~~z\in V\;.\eqno(4.3)$$
Clearly $v^{-1}(-\infty)=P$.  To show that $v$ is quasi-psh it suffices to
verify that if $z_0\in \partial \Omega_k$, then for $z\in\Omega_k\ssm P$
sufficiently close to $z_0$ we have $v_k(z)<v(z)$.  This is obvious if
$z_0\not\in P$ since $v_k\to -\infty$ as $z\to z_0$.
So consider $z_0\in \partial \Omega_k\cap P\cap \Omega'_j$.  Suppose that
$z\in \Omega_k\cap\Omega'_j\ssm P$ is sufficiently close to $z_0$ so that
$$\max\{u_j(z),u_k(z)\}<\min\{t_j,t_k\}~~~\hbox{and}~~~\varphi_k(z)\leq
-3\;.$$
Then we have $u_j(z)\in[t_{\nu+1},t_\nu[$ and
$u_k(z)\in[t_{\mu+1},t_\mu[$ with indices $\mu,\nu\geq\max\{j,k\}$.
Then by (4.2), $u_j(x)\geq t_{\mu+2}$ and $u_k(x)\geq t_{\nu+2}$; hence
$|\mu-\nu|\leq 1$ and $|\chi\circ u_j(z)-\chi\circ u_k(z)|\leq 2$.
Therefore
$$v(z)\geq v_j(z)=\chi\circ u_j(z)\geq\chi\circ u_k(z)-2\geq v_k(z)+1\;,$$
which completes the proof of (ii).

To obtain the global case of property (i), we replace the max function
in (4.3) by the regularized max functions

$${\textstyle\max_\varepsilon}=\max{}\star\rho^m_\varepsilon:
\rb^m\rightarrow \rb\;,$$
where $\rho^m_\varepsilon$ is of the form
$\rho_\varepsilon^m(x_1,\ldots,x_m)=\varepsilon^{-m}\rho(x_1/\varepsilon)
\cdots\rho(x_m/\varepsilon)$. This form of smoothing kernel ensures that
$$x_m\leq\max\{x_1,\ldots,x_{m-1}\}-2\varepsilon\Rightarrow
{\textstyle \max_\varepsilon}(x_1,\ldots,x_m)={\textstyle \max_\varepsilon}
(x_1,\ldots,x_{m-1})\,.$$
(Of course, $\max_\varepsilon $ is convex and is invariant under
permutations of variables.) The function $v_\varepsilon $ obtained by
using $\max_\varepsilon $ in (4.3) with any $\varepsilon \leq 1/2$ is
quasi-psh on $Y$, smooth on $Y\ssm P$, and has polar set $P$.
\finpr\vskip7pt

Now we show that for any Stein submanifold $M$ of a complex manifold~$Y$,
there are complete pluripolar graphs of sections $K\to N_M$ over arbitrary
compact subsets~$K\subset M$.

\than Proposition 4.6|Let $M$ be a $($not necessarily closed$)$ Stein
submanifold of a complex manifold~$Y$, with $\dim M<\dim Y$, and let
$K$ be a holomorphically convex compact subset of $M$.
Let $\psi:U\to V$ be as in Lemma~2.1{\rm(i)}.
Then there is a continuous section $g:K\to N_M$ which is
holomorphic in the interior~$K^\circ$ of~$K$, such that
$g(K)\subset U$ and the sets $K_\varepsilon= \psi(\varepsilon g(K))$ are
complete pluripolar in~$Y$ for $0<\varepsilon\leq 1$.\fintha

Here $\varepsilon g(K)$ denotes the $\varepsilon$-homothety of the
compact section $g(K)\subset N_M$; thus $K_\varepsilon$ converges to $K$ as
$\varepsilon$ tends to~$0$.

\proof By 4.5(ii), the final assertion is equivalent to proving
that $\varepsilon g(K)$ is complete pluripolar in $N_M$,
and we may of course assume that $\varepsilon =1$.  The
holomorphic convexity of $K$ implies the existence of a
sequence $(f_j)_{j\geq 1}$ of holomorphic functions on $M$
such that $K=\cap \{|f_j|\leq 1\}$.  By Lemma 2.5, there is
a finite sequence $(g_k)_{1\leq k\leq N}$ of holomorphic
sections of $N_M$ without common zeroes.  We consider the
sequence $(j_\nu, k_\nu)_{\nu\geq 1}$ of pairs of positive
integers obtained as the concatenation $\Lambda_1,
\Lambda_2,\ldots,\Lambda_\ell,\ldots$ of the finite sequences
$$\Lambda_\ell = (j_q, k_q)_{\nu(\ell)\leq
q<\nu(\ell+1)}=\left((1,1),\ldots,(1,N),\ldots,(\ell,
1),\ldots,(\ell,N)\right)$$
(where $\nu(\ell) = N\ell(\ell-1)/2+1)$.  We define $g$ to
be the generalized ``gap sequence"
$$g(z) = \sum^\infty_{\nu=1} e^{-p_{\nu}}
f_{j_{\nu}}(z)^{p_{\nu}^{2}}g_{k_{\nu}}(z)~~~~\forall z\in K\;,\eqno(4.4)$$
where $(p_\nu)$ is a strictly increasing sequence of
positive integers to be defined later.  The series (4.4)
converges uniformly on $K$ (with respect to a Hermitian
metric on $N_M$) and thus $g$ is continuous on $K$ and
holomorphic on $K^\circ$.  To construct our psh function with
polar set $g(K)$ we first select a smooth Hermitian metric
$\|~~\|$ on $N_M$ such that $\zeta \mapsto \log\|\zeta\|$
is psh on $N_M$ (e.g., take $\|\zeta\|^2 = \sum
|g^\star_i(z)\cdot \zeta|^2$ for $\zeta\in N_{M,z}$, where
$(g^\star_i)$ is a collection of holomorphic sections of $N_M^\star$
without common zeroes).  We consider the global holomorphic sections
$$s_q = \sum^q_{\nu=1} e^{-p_{\nu}} f_{j_{\nu}}^{p_{\nu}^{2}}
g_{k_\nu}\eqno(4.5)$$ of $N_M$, which converge to $g$ uniformly on
$K$ exponentially fast.  We consider the psh functions $u_\ell$ on $N_M$
given by
$$u_\ell (\zeta) = \max_{\nu(\ell)\leq q<\nu (\ell+1)}
{1\over p^4_q} \log \|\zeta - s_q (z) \|\;,~~~~\forall
\zeta\in N_{M,z},\;\; \forall z\in M\;.\eqno(4.6)$$
We shall show that the infinite sum
$$u = \sum^\infty_{\ell=1} \max (u_\ell,-1)\eqno(4.7)$$
is psh on $N_M$ with $u^{-1} (-\infty) = g(K)$.  For
$z\in K$, $\zeta = g(z) \in g(K)$, we obtain from (4.4) the
estimate
$$\|\zeta-s_q(z)\|\leq ce^{-p_{q+1}}$$
where $c$ is independent of the choice of $z\in K$.  Thus
if we choose the $p_q$ with $p_{q+1}\geq 2p^4_q$, we have
$u_\ell < -1$ on $g(K)$ for $\ell$ large, and thus
$u\equiv -\infty$ on $g(K)$.  Next we show that $u$ is
psh.  For this it suffices to show that
$\sum^\infty_{\ell=1} \sup_Au_\ell<+\infty$ for all compact
subsets $A$ of $N_M$.  Let $M_q$ be an exhausting sequence
of compact subsets of $M$.  When $\zeta\in N_{M,z}$ with
$z\in M_q$, we have
$$\log\|\zeta-s_q(z)\|\leq \log (1+\|\zeta\|) + (p^2_q + 1)
\log (1+C_q)$$
where $C_q$ is the maximum of $\|g_1\|,\ldots,\|g_N\|,
|f_1|,\ldots,|f_q|$ on $M_q$.  In addition to our previous
requirement, we assume that $p_q \geq \log (1+C_q)$ and hence
$$u_\ell (\zeta)\leq {2 + \log
(1+\|\zeta\|) \over p_{\nu(\ell)}},\;\;\forall \zeta\in
N_{M, z},\;\;\forall z\in M_{\nu(\ell)}\;.$$
Since $p_{\nu(\ell)} \geq p_\ell \geq 2^{\ell-1}$ (in fact, the growth is at
least doubly exponential), it follows that the series $\sum u_\ell$ converges
to a psh function $u$ on $N_M$.

We finally show that $u(\zeta)$ is finite outside $g(K)$. Fix a point
$\zeta \in N_{M, z} - g(K)$.  If $z\in K$, then the logarithms
appearing in (4.6) converge to $\log\|\zeta-g(z)\|>-\infty$, and hence
$u(\zeta)>-\infty$. Now suppose $z\not\in K$. Then one of the values
$f_j(z)$ has modulus greater than~$1$. Thus for $\ell$ sufficiently
large, there is a $q=q(\ell)$ in the range $\nu(\ell)\leq
q<\nu(\ell+1)$ such that $|f_{j_{q}}(z)|>1$ and $\|g_{k_{q}}(z)\| =
\max\|g_i(z)\|>0$. Assume further that $|f_{j_{q}}(z)|=\max_{1\leq
i\leq\ell}|f_i(z)|$. For this choice of $q$ we have
$$\|\zeta-s_q(z)\|\geq e^{-p_{q}}|f_{j_{q}}(z)|^{p_{q}^{2}}|g_{k_{q}}(z)|-
{|f_{j_{q}}(z)|^{p_{q-1}^2+1}|g_{k_{q}}(z)|\over|f_{j_{q}}(z)|-1}-
\|\zeta\|\;.$$
Since $p_q\geq 2p_{q-1}^4$, we easily conclude that $\|\zeta -
s_{q(\ell)} (z) \|>1$ for $\ell$ sufficiently large and hence $u_\ell
(z)$ is eventually positive. Therefore $u(z)>-\infty$; hence $u^{-1}
(-\infty)=g(K)$, and so $g(K)$ is complete pluripolar in $N_M$.\finpr

\rem Remark|If one chooses $p_\nu$ to be larger than the norm of the
first $\nu$ derivatives of $f_1,\ldots,f_\nu$ on $K$ in the above
proof, the section $g$ will be smooth on $K$.

\than Proposition 4.7|Let $K\i M\i Y$ be as in
Proposition~4.6, and let $I$ be a finite subset of the
interior $K^\circ$ of $K$.  Then for every positive integer
$m$ we can select the section $g:K\rightarrow N_M$ in
Proposition~4.6 with the additional property that the
$m$-jet of $g$ vanishes on $I$ $($and thus $K_\varepsilon $
is tangent to $K$ of order $m$ at all points of $I)$.\fintha

\proof The proof of Proposition~4.6 used only the fact that
the sections $g_1\ld g_N$ do not vanish simultaneously on
$M\ssm K$. Hence, we can select the $g_j$ to vanish at the
prescribed order $m$ at all points of $I$.\finpr\vskip7pt

Next, we prove the following simple approximation theorem based on
H\"ormander's $L^2$ estimates (see Andreotti-Vesentini [AV65] and
[H\"o66]).

\than Proposition 4.8|Let $X$ be a projective algebraic
manifold and let $L$ be an ample line bundle on~$X$. Let $P$ be a
complete pluripolar set in $X$ such that $L$ is trivial on a
neighborhood $\Omega$ of~$P$; i.e., there is a nonvanishing section
$s\in H^0(\Omega,L)$. Then there is a smaller neighborhood $V\subset
\Omega$ such that every holomorphic function $h$ on  $\Omega$ can be
approximated uniformly on $V$ by a sequence of global sections $h_\nu\in
H^0(X,L^{\otimes\nu})$; precisely, $h_\nu s^{-\nu}\to h$ uniformly
on~$V$. \fintha

\proof By Lemma 4.5(i), we can choose a quasi-psh
function $u$ on $X$ with $u^{-1}(-\infty)=P$ and $u$~smooth on $X\ssm P$.
Fix a K\"ahler metric $\omega$ on $X$ and a Hermitian metric on $L$
with positive definite curvature form $\Theta(L)=\sqrt{-1}
\partial\op\varphi$, where $e^{-\varphi}$ is the weight representing the
metric of~$L$. Precisely, let $(U_\alpha)_{0\leq\alpha\leq N}$ be an open
cover of $X$ with trivializing sections $s_\alpha\in H^0(U_\alpha,L)$.  We let
$\|s_\alpha\|^2 = e^{-\varphi_\alpha}$ so that for $v=v^\alpha s_\alpha(z)\in
L_z$ we have $\|v\|^2=|v^\alpha|^2 e^{-\varphi_\alpha(z)}$.  In the following,
we shall take $U_0=\Omega$ and $s_0=s$.  We shall also drop the index $\alpha$
and denote the square of the norm of $v$ simply by $\|v\|^2=|vs^{-1}|^2
e^{-\varphi(z)}$.

By multiplying $u$ with a small positive constant,
we may assume that $\Theta(L)+\sqrt{-1}\partial\op
u\geq\varepsilon\omega$ for some $\varepsilon>0$. Then, for $\nu\in\nb$ large,
H\"ormander's $L^2$ existence theorem (as given in [De92, Theorem 3.1])
implies that for every $\op$-closed $(0,1)$-form $g$ with values in
$L^{\otimes\nu}$ such that $$\int_X\|g\|^2e^{-\nu u}dV_\omega<+\infty$$
(with $dV_\omega={}$K\"ahler volume element), the equation $\op f=g$ admits
a solution $f$ satisfying the estimate
$$\int_X|fs^{-\nu}|^2e^{-\nu(\varphi+u)}dV_\omega
=\int_X\|f\|^2e^{-\nu u}dV_\omega\le{1\over\lambda_\nu}
\int_X\|g\|^2e^{-\nu u}dV_\omega,$$
where $\lambda_\nu$ is the minimum of the eigenvalues of
$$\nu\big(\Theta(L)+\sqrt{-1}\partial\op u\big)+{\rm Ricci}(\omega)$$
with respect to $\omega$ throughout~$X$. Clearly
$\lambda_\nu\geq\varepsilon\nu+\rho_0$
where $\rho_0\in\rb$ is the minimum of the
eigenvalues of Ricci$(\omega)$. We apply the existence theorem to the
$(0,1)$-form $g=\op(h\theta s^\nu)=h\op\theta s^\nu$, where $\theta$
is a smooth cut-off function such that
$$\eqalign{
&\theta=1~~~~\hbox{\rm on}~~\big\{z\in\Omega\,:\,\varphi(z)+u(z)\le
c\big\},\cr
&\Supp\,\theta\subset\big\{z\in\Omega\,:\,\varphi(z)+u(z)<c+1\big\}
\compact \Omega\cr}$$
for some sufficiently negative $c<0$. Then we get a
solution of the equation $\op f_\nu=h\op\theta s^\nu$ on $X$ with
$$\int_X|f_\nu s^{-\nu}|^2e^{-\nu(\varphi+u)}dV_\omega
\le{1\over\varepsilon\nu+\rho_0}
\int_X|h|^2|\op\theta|^2e^{-\nu(\varphi+u)}dV_\omega.$$
Since $\Supp(\op\theta)\subset \{c\leq\varphi(z)+u(z)\leq c+1\}$, we infer
$$e^{-\nu(c-1)}\int_{\{\varphi+u<c-1\}}|f_\nu s^{-\nu}|^2dV_\omega
\le{e^{-\nu c}\over\varepsilon\nu+\rho_0}\int_\Omega
|h|^2|\op\theta|^2dV_\omega.$$
Thus $f_\nu s^{-\nu}\to 0$ in $L^2(\{\varphi+u<c-1\})$.  Since the
$f_\nu s^{-\nu}$ are holomorphic functions on $\{\varphi+u<c-1\}$, it
follows that $f_\nu s^{-\nu}\to 0$ uniformly on the neighborhood
$V=\{\varphi+u<c-2\}\compact\Omega$ of~$P$.  Hence $h_\nu=h\theta
s^\nu-f_\nu$ is a holomorphic section of $L^{\otimes\nu}$ such that
$h_\nu s^{-\nu}$ converges uniformly to~$h$ on~$V$.\finpr

\than Corollary 4.9|Let $X$ be a projective algebraic manifold and let
$M$ be a Stein submanifold with $\dim M<\dim X$, such that there exists
an ample line bundle $L$ on $X$ with trivial restriction~$L_{|M}$.
Let $K$ be a holomorphically convex compact subset of $M$, and let
$(K_\varepsilon)_{0<\varepsilon\leq 1}$ be the complete pluripolar
sets constructed in Proposition~$4.6$. Then
for each sufficiently small $\varepsilon$, there is a positive integer
$\nu(\varepsilon)$ and a section $h_{\nu(\varepsilon)}\in
H^0(X,L^{\otimes\nu(\varepsilon)})$ such that $K_\varepsilon$ does not
intersect the divisor of~$h_{\nu(\varepsilon)}$; i.e., $K_\varepsilon$ is
contained in the affine Zariski open set
$X\ssm h_{\nu(\varepsilon)}^{-1}(0)$.\fintha

\proof By Lemma~2.4, $L$ is trivial on a neighborhood $\Omega$ of $M$.
For $\varepsilon>0$ small, we have $K_\varepsilon\subset\Omega$ and we can
apply Proposition~4.8 to get a uniform approximation $h_{\nu(\varepsilon)}
\in H^0(X,L^{\otimes\nu(\varepsilon)})$ of the function $h=1$ such that
$|h_{\nu(\varepsilon)}s^{-\nu(\varepsilon)}-h|<1/2$ near~$K_\varepsilon$.
Then $K_\varepsilon\cap
h_{\nu(\varepsilon)}^{-1}(0)=\emptyset$.\finpr\vskip7pt

We now use Corollary~4.9 to prove Theorem~4.1:

\noindent{\sl Proof of Theorem~4.1}\pointir Choose a domain $\Omega_1$ with
$\Omega_0\i\i\Omega_1\i\i\Omega$ such that $\Omega_1$ is a Runge domain in
$\Omega$ (and hence in $S$). By Corollary~4.9 applied to the compact set
$$K=\hbox{holomorphic hull of $f(\ov\Omega_1)$ in the Stein variety
$f(\Omega)\subset X$},$$ we get an embedding $j_\varepsilon: K\to X$ of
the form $j_\varepsilon(z)= \psi(\varepsilon g(z))$, converging to
$\Id_K$ as $\varepsilon$ tends to~$0$, such that $j_\varepsilon(K)$
is contained in an affine Zariski open subset~$Y_{\nu,\varepsilon}=X\ssm
h_{\nu,\varepsilon}^{-1}(0)$, $h_{\nu,\varepsilon}\in
H^0(X,L^{\otimes\mu(\nu,\varepsilon)})$. Moreover, by Proposition~4.7,
we can take $j_\varepsilon$ tangent to $\Id_{K}$ at order $k$ at
all points $t_j$. We first approximate $f$ by $F_\nu=j_{\varepsilon
_\nu}\circ f$ on $\Omega_1$, choosing $\varepsilon_\nu$ so small
that $\sup_{\Omega_1}\delta(F_\nu(z),f(z))<1/\nu$ with respect to
some fixed distance $\delta$ on~$X$. We denote
$Y_\nu=Y_{\nu,\varepsilon_\nu}$ in the sequel. By construction,
$F_\nu(\Omega_1)\subset j_{\varepsilon_\nu}(K)\subset Y_\nu$.

Let $Y_\nu\subset\cb^N$ be an algebraic embedding of $Y_\nu$ in some affine
space. Then $F_\nu$ can be viewed as a map
$$F_\nu=(F_{1,\nu}\ld F_{N,\nu}):\Omega_1\to\cb^N,~~~~
F_\nu(\Omega_1)\subset Y_\nu.$$
By Lemma~2.1, there is a retraction $\rho_\nu:V_\nu\to Y_\nu$ defined on a
neighborhood $V_\nu$ of~$Y_\nu$, such that the graph $\Gamma_{\rho_\nu}$
is contained in an algebraic variety $A_\nu$ with $\dim A_\nu=N$.
Now, since $\Omega_1$ is a Runge domain in~$S$, there are polynomial
functions $P_{k,\nu}$ in the structure ring $\cb[S]$, such that
$F_{k,\nu}-P_{k,\nu}$ is uniformly small on $\ov\Omega_0$
for each $k=1\ld N$. Of course, we can make a further finite interpolation
to ensure that $P_{k,\nu}-F_k$ vanishes at the prescribed order $m$
at all points~$t_j$. We set
$$P_\nu=(P_{1,\nu}\ld P_{N,\nu}):S\to\cb^N,~~~~f_\nu=\rho_\nu\circ P_\nu
:\Omega_0\to Y_\nu\;,$$ and take our approximations $P_{k,\nu}$
sufficiently close to $F_{k,\nu}$, in such a way that $f_\nu$ is well defined
on $\Omega_0$ (i.e. $P_\nu(\Omega_0)\subset V_\nu$), is an embedding of
$\Omega_0$ into $Y_\nu$, and satisfies $\delta(f_\nu,F_\nu)<1/\nu$
on~$\Omega_0$. Then $\delta(f_\nu,f)<2/\nu$ on~$\Omega_0$. Now,
the graph of~$f_\nu$ is just the submanifold
$$\Gamma_{f_\nu}=(\Omega_0\times Y_\nu)\cap(P_\nu\times\Id_{Y_\nu})^{-1}
(\Gamma_{\rho_\nu}).$$
Each connected component of $\Gamma_{f_\nu}$ is contained in an irreducible
component $G_{\nu,j}$ of dimension $n=\dim S$ of the algebraic variety
$(P_\nu\times\Id_{Y_\nu})^{-1}(A_\nu)$, hence $\Gamma_{f_\nu}$ is
contained in the $n$-dimensional algebraic variety $G_\nu=\bigcup
G_{\nu,j}$.\finpr

\rem Remark 4.10|The existence of an ample line bundle on $X$ with
trivial restriction to $M$ is automatic if
$H^1(M,\oc_M^\star)=H^2(M,\zb)=0$, since in that case every holomorphic
line bundle on $M$ is trivial. In Corollary~4.9, the conclusion that
$K$ admits small perturbations $K_\varepsilon$ contained in affine
Zariski open sets does not hold any longer if the hypothesis on the
existence of an ample line bundle $L$ with $L_{|M}$ trivial is omitted.
In fact, let $M$ be a Stein manifold and let $K$ be a holomorphically
convex compact set such that the image of the restriction map
$H^2(M,\zb)\to H^2(K,\zb)$ contains a nontorsion element. Then there
exists a holomorphic line bundle $E$ on $M$ such that all positive
powers $E_{|K}^{\otimes m}$ are topologically nontrivial. Let
$g_0,\dots,g_N\in H^0(M,E)$ be a collection of sections which generate
the fibres of $E$ and separate the points of $M$. By the proof of
Corollary~2.6, we find a holomorphic embedding $\varphi:M\to\pb^N$ into
some projective space, such that $E\simeq\varphi^\star\oc(1)$. Since
$\varphi^\star\oc(m)_{|K}\simeq E_{|K}^{\otimes m}$, we infer that
$\oc(m)_{|\varphi(K)}$ is topologically nontrivial for all~$m>0$. It
follows that $\varphi(K)$ cannot be contained in the complement
$\pb^N\ssm H$ of a hypersurface of degree~$m$, because $\oc(m)$ is
trivial on these affine open sets. The same conclusion holds for
arbitrary homotopic deformations $\varphi(K)_\varepsilon$
of~$\varphi(K)$.

To give a specific example, let $M = \cb^\star \times \cb^\star$ and
$K=\partial \Delta \times \partial \Delta$.  We consider the torus
$T=\cb /\{1,\sqrt{-1}\}\zb$, which has a holomorphic embedding
$g:T\rightarrow \pb^2$.  Let $\varphi:M\rightarrow \pb^8$ be the
holomorphic embedding given by $\varphi (z, w) = \psi (g(t), (1:z:w))$
where $\psi:\pb^2 \times \pb^2\rightarrow \pb^8$ is the Segre embedding
and $t = {\sqrt{-1} \over 2\pi} \log z + {1\over 2\pi} \log w$. Since
$\varphi (K)$ is a $C^\infty$ deformation of the elliptic curve $A
=\psi (g(T)\times \{(1:0:0)\})\i \pb^8$, then for any deformation
$\varphi(K)_\varepsilon $ of $\varphi(K)$ and for any algebraic
hypersurface $H\i\pb^8$ we have $\varphi(K)_\varepsilon \cdot H=A\cdot H>0$
(where $\cdot$ denotes the intersection product in the homology of
$\pb^8$), and therefore $\varphi (K)_\varepsilon \cap H\neq \emptyset$.
\vskip12pt

\bigtitle 5. Exhaustion of Stein manifolds by Runge|
domains of affine algebraic manifolds|

The methods developed in this paper are used in this section to give a
new proof of Stout's Theorem~1.6 by methods substantially different from
those in Stout [St84] and in Tancredi-Tognoli [TT89]. These methods are then
used to obtain Theorem~1.8 from Proposition~3.2. To prove Theorem~1.6, we
first show the existence of proper embeddings of arbitrarily large Stein
domains into affine algebraic manifolds, such that the normal bundle of the
embedding is trivial. Secondly, when such an embedding is given, the embedded
manifold is a global complete intersection, so it can be approximated through
an approximation of its defining equations by polynomials.

We begin with an elementary lemma about Runge neighborhoods.

\than Lemma 5.1|Let $S$ be a closed submanifold of a Stein manifold~$Y$.
Then there is a fundamental system of neighborhoods $\Omega$ of $S$ which
are Runge open sets in~$Y$.
\fintha

\proof Let $f_1=\ldots=f_N=0$ be a finite system of defining equations
of $S$ in~$Y$, and let $u$ be a strictly plurisubharmonic smooth
exhaustion function of~$Y$. For every convex increasing function
$\chi\in C^\infty(\rb)$, we set
$$\Omega_\chi=\Big\{\sum_{1\leq j\leq N}|f_j(z)|^2e^{\chi(u(z))}<1\Big\}.$$
Since $\Omega_\chi$ is a sublevel set of a global psh function on~$Y$,
we infer that $\Omega_\chi$ is a Runge domain in~$Y$ by
[H\"o66,~Th.~5.2.8]. Clearly $\{\Omega_\chi\}$ is a fundamental system of
neighborhoods of~$S$.\finpr

\than Lemma 5.2|Let $S$ be a Stein manifold. For any relatively compact
Stein domain $S_0\compact S$, there is a proper embedding $j:S_0\to A$
into an affine algebraic manifold, such that the image $j(S_0)$ has
trivial normal bundle $N_{j(S_0)}$ in~$A$.
\fintha

\proof By the Bishop-Narasimhan embedding theorem, we can suppose that
$S$ is a closed submanifold of an affine space~$\cb^p$. Since any exact
sequence of vector bundles on a Stein manifold splits, $T_S\oplus
N_{S/\cb^p}$ is isomorphic to the trivial bundle $T_{\cb^p|S}=S\times
\cb^p$. Hence, to make the normal bundle of $S$ become trivial, we need
only to embed a Stein neighborhood of $S_0$ (in $\cb^p$) into an
algebraic manifold so that its normal bundle when restricted to $S_0$
is isomorphic to~$T_S$.

For this, we choose by Lemma~2.1 a holomorphic retraction
$\rho:\Omega\to S$ from a neighborhood $\Omega$ of $S$ in $\cb^p$
onto~$S$. We can suppose $\Omega$ to be a Runge open subset of~$\cb^p$,
by shrinking $\Omega$ if necessary (apply Lemma~5.1). Let
$\Omega_0\compact\Omega$ be a Runge domain in $\cb^p$ containing~$\ov
S_0$. By Proposition~3.2 applied to the holomorphic vector bundle
$\rho^\star T_S$ on~$\Omega$, there is an open embedding $\varphi:\Omega_0\to
Z$ into a projective algebraic manifold~$Z$, an algebraic vector bundle
$E\to Z$ such that $\varphi^\star E\simeq(\rho^\star T_S)_{|\Omega_0}$,
and an ample line bundle $L\to Z$ such that $\varphi^\star L$ is
trivial. We consider the composition of embeddings
$$S\cap\Omega_0\lhra\Omega_0\buildo{\displaystyle\varphi}\over\lhra Z
\lhra E\lhra \wh E,$$
where $Z\hookrightarrow E$ is the zero section and $\wh E=P(E\oplus\cb)$
is the compactification of $E$ by the hyperplane at infinity in each fibre.
As $N_{Z/\wh E}\simeq E$, we find
$$\eqalign{
(N_{Z/\wh E})_{|S\cap\Omega_0}&\simeq(\varphi^\star E)_{|S\cap\Omega_0}
\simeq(\rho^\star T_S)_{|S\cap\Omega_0}=(T_S)_{|S\cap\Omega_0}\,,\cr
N_{S\cap\Omega_0/\wh E}&=(N_{S/\cb^p})_{|S\cap\Omega_0}\oplus
(N_{Z/\wh E})_{|S\cap\Omega_0}\cr
&\simeq(N_{S/\cb^p}\oplus T_S)_{|S\cap\Omega_0}
\simeq(S\cap\Omega_0)\times\cb^p\,.\cr}$$
Now, $\wh E$ is equipped with a canonical line bundle $\oc_{\wh E}(1)$,
which is relatively ample with respect to the fibres of the projection
$\pi:\wh E\to Z$. It follows that there is an integer $m\gg 0$ such
that $\wh L:=\pi^\star L^{\otimes m}\otimes\oc_{\wh E}(1)$ is ample
on~$\wh E$. Now, the zero section of $E$ embeds as
$P(0\oplus\cb)\subset P(E\oplus\cb) =\wh E$, hence $\oc_{\wh
E}(1)_{|Z}\simeq Z\times\cb$. Since $\varphi^\star L$ is trivial, we
infer that $\wh L_{|S\cap\Omega_0}$ is also trivial. Apply Corollary~4.9
to $X=\wh E$ with the ample line bundle $\wh L\to\wh E$, and to some
holomorphically convex compact set $K$ in the submanifold
$S\cap\Omega_0\subset\wh E$, such that~$K\supset\ov S_0$. We infer the
existence of a holomorphic embedding $i:S_0\to Y$ into an affine
Zariski open subset $Y$ of~$\wh E$. By construction, this embedding is
obtained from the graph of a section $K\to N_{S\cap\Omega_0/\wh E}\,$,
hence $$N_{S_0/Y}\simeq (N_{S\cap\Omega_0/\wh E})_{|S_0}\simeq
S_0\times\cb^p.$$ It remains to show that we can find a {\sl proper}
embedding $j$ satisfying this property. We simply set
$$j=(i,i'):S_0\la A=Y\times\cb^q,$$
where $i':S_0\to\cb^q$ is a proper embedding of $S_0$ into an affine space.
Clearly, $j$ is proper. Moreover
$$N_{j(S_0)/Y\times\cb^q}\simeq(T_Y\oplus\cb^q)/\Im(di\oplus di')$$
with $di$ being injective, hence we get an exact sequence
$$0\la\cb^q\la N_{j(S_0)/Y\times\cb^q}\la N_{i(S_0)/Y}\la 0\,.$$
As $N_{i(S_0)/Y}$ is trivial, we conclude that
$N_{j(S_0)/Y\times\cb^q}$ is also trivial.\finpr

\than Lemma~5.3|Let $S_0$ be a closed Stein submanifold of an affine
algebraic manifold $A$, such that $N_{S_0/A}$ is trivial. For any Runge
domain $S_1\compact S_0$ and any
neighborhood $V$ of $S_0$ in $A$, there exist a Runge domain
$\Omega$ in $A$ with $S_0\subset\Omega\subset V$, a holomorphic retraction
$\rho:\Omega\to S_0$, and a closed algebraic submanifold $Y\subset A$ with
the following properties: \smallskip
\item{\rm(i)} $Y\cap\rho^{-1}(S_1)$ is a Runge open set in
$Y$;\vskip0pt
\item{\rm(ii)} $\rho$ maps $Y\cap\rho^{-1}(S_1)$
biholomorphically onto $S_1$. \fintha

\proof By Lemmas~2.1 and 5.1, we can find a Runge domain $\Omega$ in $A$
with $S_0\subset\Omega\subset V$ and a holomorphic retraction
$\rho:\Omega\to S_0$. Let $\ic$ be the ideal sheaf of $S_0$ in~$A$. Then
$\ic/\ic^2$ is a coherent sheaf supported on~$S_0$, and its restriction to
$S_0$ is isomorphic to the conormal bundle $N_{S_0/A}^\star$. Since this
bundle is trivial by assumption, we can find $m=\codim S_0$ global generators
$g_1\ld g_m$ of $\ic/\ic^2$, as well as liftings $\wt g_1\ld\wt g_m\in
H^0(A,\ic)$ ($A$ is Stein, thus $H^1(A,\ic^2)=0$). The system of equations
$\wt g_1 =\ldots=\wt g_m=0$ is a regular system of equations for $S_0$ in a
sufficiently small neighborhood. After shrinking~$\Omega$, we can suppose that
$$S_0=\big\{z\in\ov\Omega\,:\,\wt g_1(z)=\ldots=\wt g_m(z)=0\big\},~~~~
d\wt g_1\wedge\ldots\wedge d\wt g_m\ne 0~~\hbox{\rm on}~\ov\Omega.$$
Let $P_{j,\nu}$ be a sequence of polynomial functions on~$A$ converging
to $\wt g_j$ uniformly on every compact subset of~$A$. By Sard's
theorem, we can suppose that the algebraic varieties $Y_\nu=\bigcap_
{1\leq j\leq m}P_{j,\nu}^{-1}(0)$ are nonsingular (otherwise, this can
be obtained by adding small generic constants $\varepsilon_{j,\nu}$
to~$P_{j,\nu}$). It is then clear that $Y_\nu\cap\Omega$ converges to
$S_0$ uniformly on all compact sets. In particular, the restriction
$$\rho_\nu:Y_\nu\cap\rho^{-1}(S_1)\la S_1$$
of $\rho$ is a biholomorphism of $Y_\nu\cap\rho^{-1}(S_1)$ onto $S_1$
for $\nu$ large. This shows that $S_1$ can be embedded as an open subset
of $Y_\nu$ via~$\rho_\nu^{-1}$.

It remains to show that $Y_\nu\cap\rho^{-1}(S_1)$ is a Runge domain
in~$Y_\nu$ for $\nu$~large. Let $h$ be a holomorphic function
on~$Y_\nu\cap\rho^{-1}(S_1)$. Then $h\circ\rho_\nu^{-1}$ is a
holomorphic function on~$S_1$. Since $S_1$ is supposed to be a Runge
domain in~$S_0$, we have $h\circ\rho_\nu^{-1}=\lim h'_k$ where $h'_k$
are holomorphic functions on $S_0$. Then $h'_k\circ\rho$ is holomorphic
on~$\Omega$, and its restriction to $Y_\nu\cap\rho^{-1}(S_1)$ converges
to $h$ uniformly on compact subsets. However, $\Omega$ is also a Runge
domain in~$A$, so we can find a holomorphic (or even algebraic)
function $h_k$ on~$A$ such that $|h_k-h'_k\circ\rho|<1/k$ on
$\ov{Y_\nu\cap\rho^{-1}(S_1)}$, which is compact in $\Omega$ for $\nu$~large.
Then the restriction of $h_k$ to $Y_\nu$ converges to $h$ uniformly on all
compact subsets of~$Y_\nu\cap\rho^{-1}(S_1)$, as desired.\finpr\vskip7pt

\noindent{\sl Proof of Theorem~1.6}\pointir
Take a Stein domain $S_0\compact S$ such that $\Omega\compact S_0$.
Then $S_1=\Omega$ is also a Runge domain in~$S_0$, and Theorem~1.6 follows
from the combination of Lemmas~5.2 and~5.3.\finpr\vskip7pt

We can now use Theorem~1.6 and Lemmas~5.2 and~5.3 to infer
Theorem~1.8 from Proposition~3.2.

\noindent{\sl Proof of Theorem~1.8}\pointir By Theorem~1.6, we can take
a Runge domain $\wt\Omega\compact S$ such that
$\wt\Omega\supset\!\supset \Omega$ and embed $\wt\Omega$ as a Runge
domain in an affine algebraic manifold~$M$. By applying Proposition~3.2 with
$M$ playing the role of $S$, $\wt\Omega$ the role of $\Omega$, and $\Omega$
the role of $\Omega_0$, we obtain a projective algebraic
manifold $Z$ containing $\Omega$ and an algebraic vector bundle $\wt E\to
Z$ such that $\wt E_{|\Omega}\simeq E_{|\Omega}$. Moreover, there is an ample
line bundle $L$ on $Z$ such that $L_{|\Omega}$ is trivial.

Only one thing remains to be proved, namely that $\Omega$ can be taken
to be a Runge subset in an affine open subset $Z\ssm D$. (As mentioned
earlier, this would be immediate if Problem~3.5 had a positive answer.)
To reach this situation, we modify the construction as follows.
Embed $\Omega$ as $\Omega\times\{0\}\subset Z\times\pb^1$. This embedding
has trivial normal bundle and moreover
$L\stimes\oc_{\pb^1}(1)_{|\Omega\times\{0\}}$ is trivial. By
Corollary~4.9, for every holomorphically convex compact
set~$K\subset\Omega$, there is a holomorphic map
$g:K^\circ\to\cb\subset\pb^1$ which is smooth up to the boundary, such that
the graph of $g$ is contained in an affine open set $Z_1\ssm D_1$, where
$D_1$ is an ample divisor of $Z_1:=Z\times\pb^1$ in the linear system
$m|L\stimes\oc_{\pb^1}(1)|$, $m\gg 0$. Let $\Omega'\subset K^\circ$ be a
Runge open subset of~$\Omega$ and let $j_N:\Omega'\to\cb^N$ be a proper
embedding. We consider the embedding
$$j'=(\Id_{\Omega'},g,j_N):\Omega'\la Z_2:=Z\times\pb^1\times\pb^N.$$
Its image $j'(\Omega')$ is a closed submanifold of the affine open set
$(Z_1\ssm D_1)\times\cb^N=Z_2\ssm D_2$, where
$D_2=(D_1\times\pb^N)+(Z_1\times \pb^{N-1}_\infty)$ is ample in~$Z_2$.
Moreover, the normal bundle of $j$ is trivial (by the same argument as at the
end of the proof of Lemma~5.2). Let $\Omega''\compact \Omega'$ be an
arbitrary Runge domain in~$\Omega'$. By Lemma~5.3, we can
find a nonsingular algebraic subvariety $Z''\subset Z_2$ such that
$\Omega''$ is biholomorphic to a Runge open set in $Z''\ssm (Z''\cap D_2)$,
via a retraction $\rho$ from a small neighborhood of $j'(\Omega')$
onto~$j'(\Omega')$. Let $E_2=\pr_1^\star\wt E\to Z_2$ and let
$E''$ be the restriction of $E_2$ to~$Z''$. Since $\rho^\star E_2\simeq E_2$
on a neighborhood of $j'(\Omega')$ by Lemma~2.4, we
infer that $E''_{|\Omega''}\simeq E_{|\Omega''}$. Therefore the conclusions
of Theorem~1.8 hold with $\Omega''$, $Z''$, $D''=Z''\cap D_2$ and $E''$ in
place of $\Omega$, $Z$, $D$ and $\wt E$. Since $\Omega''$ can be taken to
be an arbitrary Runge open set in~$S$, Theorem~1.8 is proved.\finpr
\vskip30pt

\centerline{\bf References}
\vskip8pt\eightpoint
\lettre Oka37|
\article AV65|Andreotti, A., Vesentini, E|Carleman estimates for the
Laplace-Beltrami equation in complex manifolds|Publ.\ Math.\
I.H.E.S.|25|1965|81-130|
\article Bi61|Bishop, E|Mappings of partially analytic spaces|Amer.\ J.
Math.|83|1961|209-242|
\article CG75|Cornalba, M., Griffiths, P|Analytic cycles and vector bundles on
non-compact algebraic varieties|Invent.\ Math.|28|1975|1-106|
\article De90|Demailly, J.-P|Cohomology of $q$-convex spaces in top
degrees|Math.\ Zeit\-schrift|203|1990|283-295|
\livre De92|Demailly, J.-P|Singular Hermitian metrics on
positive line bundles|Complex algebraic varieties, Lecture
Notes in Math., Vol.\ 1507, Springer-Verlag, Berlin|
1992{\rm, 87-104}|
\divers De93|Demailly, J.-P|Algebraic criteria for the hyperbolicity of
projective manifolds|manuscript in preparation|
\livre Ei70|Eisenman, D|Intrinsic measures on complex manifolds and
holomorphic mappings|Mem.\ Amer.\ Math.\ Soc., No.~96|{\rm(}1970{\rm)}|
\article EG92|Eliashberg, Y., Gromov, M|Embeddings of Stein manifolds of
dimension $n$ into the affine space of dimension $3n/2+1$|Ann.\ of
Math.|136|1992|123-135|
\article Gr58|Grauert, H|Analytische Faserungen \"uber
holomorphvollst\"andigen R\"aumen|Math.\ Ann.|135|1958|263-273|
\livre GA74|Griffiths, P., Adams, J|Topics in algebraic and analytic
geometry|Princeton Univ.\ Press, Princeton, NJ|1974|
\livre GR65|Gunning, R.C., Rossi, H|Analytic functions of several complex
variables|Prentice-Hall, Englewood Cliffs, NJ|1965|
\article Hi64|Hironaka, H|Resolution of singularities of an algebraic variety
over a field of characteristic zero|Ann.\ of Math.|79|1964|109-326|
\livre H\"o66|H\"ormander, L|An introduction to complex analysis in several
variables|1st edition, Elsevier Science Pub., New-York, {\oldstyle 1966}; 3rd
revised edition, North-Holland Math.\ Library, Vol.\ 7, Amsterdam|1990|
\article Ko67|Kobayashi, S|Invariant distances on complex
manifolds and holomorphic mappings|J.~Math. Soc. Japan|19|1967|460-480|
\livre Ko70|Kobayashi, S|Hyperbolic manifolds and
holomorphic mappings|Marcel Dek\-ker, New York|1970|
\article La86|Lang, S|Hyperbolic and diophantine analysis|Bull.\ Amer.\
Math.\ Soc.|14|1986|159-205|
\livre La87|Lang, S|Introduction to complex hyperbolic spaces
and diophantine analysis|Springer-Verlag, Berlin|1987|
\article Na60|Narasimhan, R|Embedding of holomorphically complete
spaces|Amer.\ J. Math.|82|1960|913-934|
\article Nh52|Nash, J|Real algebraic manifolds|Ann.\ of
Math.|56|1952|405-421|
\livre NO90|Noguchi, J., Ochiai, T|Geometric function
theory|American Mathematical Society, Providence, RI|1990|
\article Oka37|Oka, K|Sur les fonctions de plusieurs variables,~II. Domaines
d'holomorphie|J.~Sci.\ Hiroshima Univ.|7|1937|115-130|
\article Ri68|Richberg, R|Stetige streng pseudokonvexe Funktionen|Math.\
Ann.|175|1968|257-286|
\livre Ro71|Royden, H|Remarks on the Kobayashi metric|Proceed.\ Maryland
Conference on Several Complex Variables, Lecture Notes in Math., Vol.~185,
Springer-Verlag, Berlin|1971|
\article Sib85|Sibony, N|Quelques probl\`emes de prolongement de courants en
analyse complexe|Duke Math.~J.|52|1985|157-197|
\article Siu76|Siu, Y.-T|Every Stein subvariety has a Stein
neighborhood|Invent.\ Math.|38|1976|89-100|
\livre St84|Stout, E.L|Algebraic domains in Stein manifolds|Proceedings
of the Conference on Banach algebras and several complex variables
(New Haven, {\oldstyle 1983}), Contemporary Math., Vol.~32, Amer.\ Math.\
Soc., Providence, RI|1984|
\article TT89|Tancredi, A., Tognoli, A|Relative approximation
theorems of Stein manifolds by Nash manifolds|Boll.\ Un.\ Mat.\ Ital.,
Serie-VII|3-A|1989|343-350|
\article TT90|Tancredi, A., Tognoli, A|On the extension of Nash
functions|Math.\ Ann.|288|1990|595-604|
\article TT93|Tancredi, A., Tognoli, A|Some remarks on the classification
of complex Nash vector bundles|Bull.\ Sci.\ Math|117|1993|173-183|
\article We35|Weil, A|L'int\'egrale de Cauchy et les fonctions de
plusieurs variables|Math.\ Ann.|111|1935|178-182|
\livre Wh65|Whitney, H|Local properties of analytic varieties|Differential
and combina\-torial topology (a symposium in honor of Marston Morse),
Princeton University Press, Princeton, NJ|1965{\rm, 205-244}|
\article Yau78|Yau, S.-T|A general Schwarz lemma for K\"ahler
manifolds|Amer.\ J. Math.|100|1978|197-203|

\tenpoint$~$\hfill (October 8, 1993)
\end